\documentclass[12pt]{article}
\usepackage{latexsym,url}
\usepackage{amsmath,amssymb}
\usepackage{graphicx}
\input{epsf.sty}	

\input xy
\xyoption{all}

\newcommand{\WriteBox}{\hfill $\Box$}	

\newcommand{\WriteBoxLN}{\WriteBox \newline}
\newcommand{\WriteTriangleLN}{\hfill $\triangle$ \newline}
\newtheorem{propositionCommand}{Proposition}
\newtheorem{theorem}[propositionCommand]{Theorem}
\newtheorem{lemma}[propositionCommand]{Lemma}

\newtheorem{corollary}[propositionCommand]{Corollary}
\newtheorem{definition}[propositionCommand]{Definition}

\newcommand{\cherries}{{\textsf C}}
\newcommand{\ignore}[1]{}
\newcommand{\N}{{\mathbb N}}
\newcommand{\textProof}{\hbox{\textbf{Proof. }}}
\newcommand{\textDiscussion}{\textbf{Discussion: }}

\newenvironment{proposition}{\begin{propositionCommand}}{\end{propositionCommand}}
\newenvironment{discussion}{\textDiscussion}{\WriteTriangleLN}
\newenvironment{proof}{\textProof}{\WriteBoxLN}

\begin{document}
\title{Probabilities on cladograms:\\ introduction to the alpha model}
\author{Daniel J. Ford}
\date{}
\maketitle
\begin{abstract}
This report introduces the alpha model.  The alpha model is a one parameter family of probability models on cladograms (binary leaf-labeled trees) which interpolates continuously between the Yule, Uniform and Comb distributions.  The single parameter $\alpha$ varies from $0$ to $1$, with $\alpha=0$ giving the Yule model, $\alpha = 1/2$ the Uniform, and $\alpha = 1$ the Comb.  For each fixed $\alpha$, the alpha model is a sequence, $\{P_n\}_{n\in \N}$, with $P_n$ a probability on cladograms with $n$ leaves.  This sequence is sampling consistent, roughly meaning that choosing a random tree from $P_n$ and deleting $k$ random leaves gives a random tree from $P_{n-k}$.  It is also Markovian self-similar.  The only other known family with these properties is the beta model of Aldous.  An explicit formula is given to calculate the probability of a given tree shape under the alpha model.
Statistics such as the expected depth of a random leaf are shown to be $O(n^\alpha)$ for $\alpha \ne 0$.
The number of cherries on a random alpha tree is shown to be asymptotically normal with known mean and variance.  Finally the shape of published phylogenies is examined, using trees from Treebase.
\end{abstract}
\tableofcontents

\section{Introduction}
This report introduces a family of probability models on cladograms, collectively called the alpha model.  Each model consists of a sequence of probabilities, one for each size of tree, which is Markovian self-similar and deletion stable.  The family of models is parameterized by a single number $\alpha\in[0,1]$ and interpolates continuously between the three most popular models on cladograms called the Yule, Uniform and Comb models.  Analogous families of models are defined for several other types of tree.

A cladogram is a rooted binary tree with $n$ leaves labeled $1$ up to $n$, a root vertex and $n-1$ internal vertices.  Cladograms  are used in biological systematics to represent the evolutionary relationship between $n$ species.  They are sometimes called phylogenetic trees, although some authors reserve this term for cladograms with edge lengths.

The three most popular probability models on cladograms are the {\em Yule model}, the {\em Uniform model} and the {\em Comb model}.  The Yule model is also referred to as the {\em neutral evolution model}.  The Uniform model assigns the uniform probability measure to cladograms of each size.
The Comb model assigns probability $1$ to the most asymmetric tree of each size.

These have the property that they are deletion stable, also called sampling consistent, and Markovian self-similar.  Informally, deletion stability means that deleting a random leaf from a random tree with $n$ leaves gives a random tree from the same model with $n-1$ leaves.  Markovian self-similarity means that the subtree below an edge is distributed independently according to the same model.  Symmetry under permutation of leaf labels is also desirable.

Previously, David Aldous has introduced a one dimensional continuous family of models, collectively called the beta model (\cite{Aldous-1996},\cite{Aldous-2000},\cite{Aldous-2001}), which interpolates between the Yule, Uniform and Comb models.  These are also deletion stable and Markovian self similar, and display qualitatively different behaviors for different values of the parameter $\beta$.

The alpha model introduced here has a very simple definition which allows many of its properties to be exactly calculated for finite values of $n$.
Basically, leaves are inserted one after another until the desired number is reached.  A leaf is inserted at a given internal edge with probability $\frac{\alpha}{n-\alpha}$ and at a given leaf edge with probability $\frac{1-\alpha}{n-\alpha}$.
Setting $\alpha=0$ results in the Yule model, $\alpha=\frac{1}{2}$ gives the Uniform model and $\alpha=1$ the Comb model.

Section $2$ introduces the necessary basic definitions and results about trees.  Four particular types of tree are defined: cladograms, fat cladograms, tree shapes and fat tree shapes.  The are related by the maps which forget leaf labels or the ordering of children.  The operation of joining two trees at the root is also defined.

In Section $3$ the alpha model is defined.  In fact, a model is defined for each of the four types of tree discussed in Section $2$.  These are related by the operations of forgetting leaf labels or ordering of children.  Markovian self-similarity and deletion stability are defined, and the alpha model shown to have these properties.  Necessary and sufficient conditions for a Markovian self-similar model to be deletion stable are also derived.  The probability of a tree under a Markovian self-similar model is calculated for each type of tree and these results applied to the alpha model.

Next, the alpha model is shown to pass through the Yule, Uniform and Comb models.  The beta model is also briefly described and shown to be different from the alpha model except where they intersect at the Yule, Uniform and Comb models.

In Section \ref{chapter:Sackins-Colless-index}, two statistics on trees are discussed.  These are Sackin's index and Colless' index.  Sackin's index is the sum of the distances from each leaf to the root, and Colless' index is the sum of the differences of the number of leaves to the left and right of each branch-point.  These are shown to differ by at most $\frac{n}{2}\log_2 n$ on a rooted binary tree with $n$ leaves.
For the alpha model, Sackin's index is shown to be $O(n^{1+\alpha})$ for $\alpha \in (0,1]$.  Thus the covariance of Sackin's index and Colless' index is asymptotically $1$ for $\alpha \in (0,1]$.  The case of the Yule model, $\alpha=0$, has been studied before.  In that case both Sackin's and Colless' index are $O(n\log n)$ with known constants and covariance.

Another statistic for cladograms or binary trees is the number of cherries, addressed in Section \ref{chapter:cherries}.  A cherry is a pair of leaves which are adjacent to each other.  McKenzie and Steel \cite{McKenzie-Steel-2000} have shown that for the Yule and Uniform models the number of cherries is asymptotically normal, with known mean and variance.  These results are extended in Section \ref{chapter:cherries} to show that, for any $\alpha\in [0,1)$, the number of cherries in a random tree from the alpha model is asymptotically normal with known mean and variance.  The Comb model, $\alpha=1$, is deterministic with exactly one cherry for a comb tree with at least 2 leaves.

Section \ref{chapter:treebase} looks at the shape of published phylogeny.
Despite the increase in published phylogeny, this appears to be the first systematic study of the shape of a large number of published phylogenetic trees, perhaps with the exception of \cite{Heard-1992}.

Natural questions to ask about the shape of cladograms or phylogentic trees include:
Are they symmetrical and flat, or asymmetrical and deep?  Is there systematic bias in reconstruction algorithms?
The trees analyzed are those in Treebase \cite{Treebase}, a free database of published phylogeny.
In the past, a major stumbling block was the lack of a measure of imbalance which could be compared across trees of different sizes, see \cite{Heard-1992} for example.  Fortunately, the maximum likelihood estimate of $\alpha$ is such a measure of imbalance.

All binary trees from Treebase (as of Nov.2004) are analyzed and their shapes compared using the alpha model.  A variety of statistics are used to consider the goodness of fit of the alpha model to this data.
Two common models for cladograms are the Yule and Uniform.  It has often been noted that published trees tend, on average, to be less balanced than Yule trees but more balanced than Uniform trees.  This observation is verified and quantified for a large set of trees.

This analysis of Treebase was carried out in November 2004 and presented at the Annual New Zealand Phylogenetics Conference in Feburary 2005, along with a brief summary of Sections 2-5.

Finally, I would like to thank my advisors Persi Diaconis and Susan Holmes who have offered much guidance and support.  This work forms part of my PhD thesis and grew out of a homework exercise in a combinatorics class of Persi's.  The analysis of Treebase was suggested by Susan Holmes.  This work was supported in part by NSF award \#0241246 (Principal investigator Susan Holmes).

\section{Basic definitions and constructions for trees}

The basic objects discussed throughout this work are trees.  These will usually have a root vertex and leaf labels. 

\begin{figure}[htbp]
    \begin{center}
    \resizebox{10cm}{!}{\includegraphics{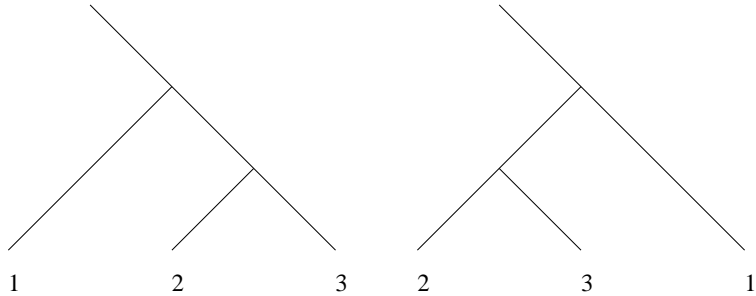}}\\
    \caption[Two representations of the same cladogram]{\label{fig:same-cladogram}The same thin cladogram, but different fat cladograms}
    \end{center}
\end{figure}

Trees will be thought of as growing down from the root.  The descendents of a vertex are those vertices further from the root, and the ancestors those which are closer to the root.  The parents and children of a vertex are those vertices immediately above and below, respectively.

Some trees are also `fat trees', in which case the ordering of the children of each vertex is important.  Thus, for fat binary trees it make sense to talk of the left and right child, and the left and right subtree below a vertex.  For thin (non-fat) trees, the children of a vertex are not ordered.  So, for example, in Figure \ref{fig:same-cladogram} the two diagrams represent the same thin tree, but different fat trees.

Isomorphisms between trees are what you might expect: graph isomorphisms which preserve any additional structure.  Isomorphic trees are considered equal.

The obvious forgetful maps which forget either labelings, or the ordering / orientation in fat trees, will also be used.

The four main type of trees considered here are:
\begin{itemize}
\item {\em tree shapes}, which are unlabeled binary rooted trees;
\item {\em cladograms}, which are tree shapes where the $n$ leaves have distinct labels $1$ up to $n$;
\item {\em fat tree shapes}, which are tree shapes where the children of each vertex are ordered; 
\item {\em fat cladograms}, which are cladograms where the children of each vertex are ordered.
\end{itemize}
The forgetful maps send each of these types to another.

The symmetric group on a labeling set acts in the obvious way on a leaf-labeled tree: by permuting the leaf labels.

Also, the {\em subtree below an edge} is defined to be the subtree consisting of all vertices and edges below and including the specified edge.  The {\em root join} of two trees is the tree formed by gluing their two root vertices together and gluing a new root edge to this vertex.  Thus the old root vertices are now the same immediate descendent of the new root vertex.

The remainder of this section is devoted to the rigorous definition of these ideas.

Finally, some familiarity with basic probability is assumed.  If you lack this background, despair not.  Most sets of interest here are finite, in which case a probability is simply a positive real function on the set which sums to $1$. {\em Independence} allows  probabilities to be multiplied in the most natural way.

Functions between finite sets extend by linearity to functions between the probabilities on these sets.  For notational convenience, the original function and its linear extension will usually be conflated.

\subsection{Graphs, trees and roots}

A {\em graph} is a pair of sets $(V,E)$, where $E \subset \{ \{u,v\} | u,v \in V\}$.  The set $V$ is called the set of {\em vertices}, and $E$ is called the set of {\em edges}.
Call $\{u,v\} \in E$ an {\em edge from $u$ to $v$}.
Note that this definition does allow 'self edges' but not 'multiple edges'.
Call $\{u,v\} \in E$ a {\em self-edge} if $u=v$.

Say that $u,v \in V$ are {\em adjacent}, or {\em neighbors}, in graph $(V,E)$ if $\{u,v\} \in E$.

A {\em path} from vertex $v_1\in V$ to vertex $v_2 \in V$ in a graph $(V,E)$ is a finite non-empty sequence $(a_i)_{i=0}^n$ such that $a_i \in V$, $a_0 = v_1$, $a_n=v_2$, and $\{a_i,a_{i+1}\} \in E$ for all $i\in \{0,1,\ldots,n-1\}$. 

The {\em length} of a path $(a_i)_{i=0}^n$ is defined to be $n$.
Note that any sequence of vertices of length $1$ is a path of length $0$.  Thus for every vertex there is a path from it to itself.

A path, $(a_i)_{i=0}^n$, is called {\em self-intersecting} if $a_i = a_j$ for some $i\ne j$.
\begin{proposition}
If there is a self-intersecting path from vertex $u$ to $v$ then there is a non self-intersecting path from $u$ to $v$.
\end{proposition}
\begin{proof}
Suppose $(a_i)_{i=0}^n$ is a self-intersecting path from $u$ to $v$.  Let
$$j = \min\left\{i | a_i = a_k, i,k\in\{0,\ldots,n\} , i\ne k \right\}$$
Choose $k\ne j$ such that $a_k = a_j$.  Now $(a_i)_{i = 1,\ldots,j,k+1,\ldots,n}$ is a path from $u$ to $v$, as $\{a_k,a_{k+1}\} \in E$ and so $\{a_j,a_{k+1}\}\in E$.  If this new path is self-intersecting then the same argument may be applied to it.  As the path length is decreased each time, this process may be repeated only a finite number of times after which the resulting path from $u$ to $v$ must be non self-intersecting. 
\end{proof}

A {\em tree} is a graph such that for each pair of vertices there is exactly one non self-intersecting path from the first vertex to the second.

The {\em distance} between two vertices is defined to be minimal length of a path from one to the other.  In other words, the distance between vertices $v_1,v_2 \in V$ is defined to be $d(v_1,v_2) = \min \{ n | (a_i)_{i=0}^n \text{ is a path from $v_1$ to $v_2$ } \}$.  Note that $\min \phi = +\infty$
\begin{proposition}
$d(\cdot,\cdot)$ is a metric.
\end{proposition}
\begin{proof}
If $(a_i)_{i=0,\ldots,n}$ is a path from $u$ to $v$ then $(a_i)_{i=n,\ldots,0}$ is a path from $v$ to $u$.  Thus $d(u,v) = d(v,u)$.  If $(a_i)_{i=0,\ldots,n}$ is a path from $u_1$ to $u_2$ of length $n$ and $(b_i)_{i=1,\ldots,m}$ is a path from $u_2$ to $u_3$ of length $m$ then $a_n=b_0$ and so $(a_0,\ldots,a_n,b_1,\ldots,b_m)$ is a path from $u_1$ to $u_3$ of length $m+n$.  Thus $d(u_1,u_3) \le d(u_1,u_2) + d(u_2,u_3)$.  Finally, $(v)$ is a path from $v$ to $v$ of length $0$ and so $d(v,v)=0$.
\end{proof}

Call a graph {\em connected} if there is a path from every vertex to every other vertex.  In other words, $d(v_1,v_2) < \infty$ for all $v_1,v_2 \in V$.

The {\em degree} of a vertex $v\in V$ in a graph $(V,E)$ is defined to be $d(v) = |\{ \{u,v\} \in E\}| + |\{ \{v,v\} \in E\}|$.  Note that a self-edge, if it exists, is counted twice.
In other words, the degree of a vertex is the number of 'half-edges' which are incident to it.

A {\em leaf} is a vertex of degree $1$.

A {\em binary tree} is a tree where every vertex has degree $1$ or $3$.  This is sometimes called a {\em trivalent tree}.

A {\em rooted graph} is a tuple $(V,E,r)$ such that $(V,E)$ is a graph and $r\in V$.  The vertex $r$ is called the {\em root} of the graph.  The empty graph $(\phi,\phi)$ may be considered as a rooted graph.  

A {\em rooted tree} is a rooted graph which is a tree, such that the root vertex is a leaf.
In this case, the set of leaves and number of leaves will {\em not} include the root vertex.  This convention will sometimes be highlighted by use of the term {\em non-root leaves}.

\subsection{Ancestors, descendents, parents and children}
Say that a path $(a_i)_{i=0}^n$ {\em passes through} vertex $x\in V$ if $a_i = x$ for some $i \in \{0,\ldots,n\}$.

For the remainder of this section, let $x$ and $y$ be vertices of a rooted tree with vertex set $V$, edge set $E$ and root $r$.

Call $y$ an {\em ancestor} of $x$ if $y$ lies on the unique non self-intersecting path from $x$ to the root.

Call $x$ a {\em descendent} of $y$ if $y$ is an ancestor of $x$.

Call $y$ a {\em parent} of $x$ if $\{u,v\}$ is an edge and $y$ is an ancestor of $x$.  Uniqueness of the non self-intersecting path from $x$ to the root and the absence of cycles implies that the parent of a vertex is unique.

Call $x$ a {\em child} of $y$ if $y$ is the parent of $x$.

In this way, the vertices of a rooted tree have a poset structure, with the root as the unique maximum element.  In this partial order, a vertex $x$ is said to be greater than a vertex $y$ if and only if $x$ is an ancestor of $y$.

Given a set of vertices, $s$, define the {\em latest common ancestor} of these vertices to be a vertex which has every element of $s$ as a descendant, but for which no descendent of this vertex has that property.  The finite tree structure guarantees that this vertex exists and is unique.

\subsection{Fat, thin, labeled, unlabeled and the forgetful maps}
The additional properties {\em fat, thin, labeled} and {\em unlabeled} are now defined, as well as the associated forgetful maps.

A {\em partial function} between two sets $X$ and $Y$ consists of a subset $Z$ of $X$ and a set map from $Z$ to $Y$.  The subset $Z$ is called the {\em domain} of the partial function.

A {\em partial labeling} of a graph is a partial function from the vertex set to a set which is called the {\em set of labels}.
A vertex is said to be {\em labeled} if it is in the domain of this partial function.

A tree together with a labeling is called a {\em labeled tree}.  If every vertex is labeled then the tree is said to be {\em totally labeled}.  Throughout this text, labelings are not assumed to be total, and partially labeled trees may be referred to simply as labeled trees.

A tree is said to be {\em leaf labeled} if it has a labeling such that the set of labeled vertices is exactly the set of leaves.  In other words, the domain of the labeling function is the set of leaves.

An {\em orientation} of a graph is a map which assigns to each vertex a cyclic ordering on its set of neighbors.  Call the image of a vertex under this map the orientation at that vertex.

A tree together with an orientation is called a {\em fat tree}, or {\em ribbon tree}.  A tree without an orientation is called a {\em thin tree}.  Trees are assumed to be thin unless stated otherwise.

The map $F_o$ forgets orientations and the map $F_l$ forgets labelings.  Thus applying $F_o$ to a partially labeled fat tree gives a partially labeled thin tree.  Applying $F_l$ to a labeled tree gives the same tree without its labeling function.
Explicitly:

\begin{definition}
\label{def:forgetfuls}
If $t$ is a fat tree then $F_o(t)$ is a thin tree with the same vertex set, edge set, and any other properties such as root or labeling.

If $t$ is a labeled tree then $F_l(t)$ is an unlabeled tree with the same vertex set, edge set, and any other properties such as root or orientation.
\end{definition}

Note that $F_o$ and $F_l$ commute, in the sense that applying $F_oF_l$ or $F_lF_o$ to a fat labeled tree gives the thin unlabeled tree with the same vertex and edge set, and any other properties such as a root.
\begin{equation}
\label{eqn:forgetfuls-commute}
F_oF_l = F_lF_o
\end{equation}

\subsection{Four types of tree: fat and thin cladograms and tree shapes}
\begin{definition}
A {\em cladogram with $n$ leaves} is a partially labeled rooted binary tree with $n$ leaves (not including the root) and label set $\{1,2,\ldots,n\}$, such that the labeled vertices are exactly the (non-root) leaves and no two leaves have the same label.  Thus each label $1,2,\ldots,n$ appears exactly once.  Define the empty labeled tree to be a cladogram with $0$ leaves.
\end{definition}

The four types of tree of particular interest here are:
\begin{itemize}
\item {\em rooted binary trees}, also called {\em tree shapes};
\item {\em cladograms}, as defined above;
\item {\em fat rooted binary trees}, also called {\em fat tree shapes}.
\item {\em fat cladograms}, which are cladograms together with an orientation.
\end{itemize}

Thus the map $F_o$ sends fat tree shapes to tree shapes, and sends fat cladograms to cladograms.  The map $F_l$ sends cladograms to tree shapes and sends fat cladograms to fat tree shapes.

\subsection{Isomorphisms of trees}
In this section, isomorphism is defined for various types of tree.  In summary, an isomorphism here is a graph isomorphism which preserves any additional structure.  General morphisms of trees are omitted, but may be easily guessed at.

An isomorphism between trees $(V_1,E_1)$ and $(V_2,E_2)$ is a bijection $f:V_1 \rightarrow V_2$ such that $\{f(u),f(v)\} \in E_2$ if and only if $\{u,v\} \in E_1$.

If either of the trees has extra structure such as a root, orientation or labeling then both must have this extra structure and it must be preserved by the map $f$.  In particular:
\begin{itemize}
\item If $r_1$ is the root of the first tree then $f(r_1)$ is the root of the second;
\item If $g_2$ is the labeling of the second tree then $fg_2$ is the labeling of the first tree;
\item If $(v_1,v_2,\ldots,v_k)$ is the cyclic orientation at vertex $v$ then \\ $(f(v_1),\ldots,f(v_k))$ is the cyclic orientation at vertex $f(v)$.
\end{itemize}

Isomorphic trees are considered equal.

\subsection{The action of the symmetric group on leaf labels}

\begin{definition}
If $t$ is a labeled tree with labeling partial function $g$ and label set $L$ and $\sigma$ is a permutation of the set $L$ then define $\sigma(t)$ to be a tree identical to $t$ except that it has labeling function $\sigma g$.
\end{definition}
In other words, apply the permutation to each label.  This defines a group action.

\begin{figure}[htbp]
    \begin{center}
    \resizebox{10cm}{!}{\includegraphics{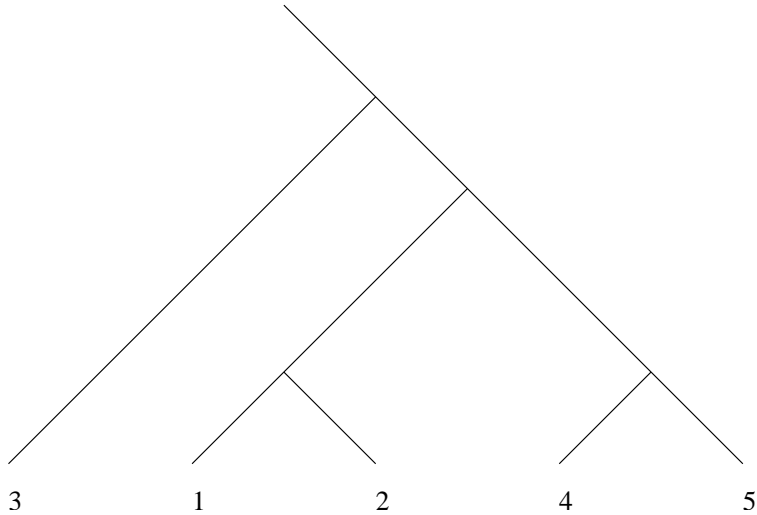}}\\
    \caption{\label{fig:perm_example}A cladogram invariant under permutations  $(12)(3)(45)$ and  $(14)(25)(3)$}
    \end{center}
\end{figure}

Some permutations will act trivially on some cladograms, such as the permutations $(12)(3)(45)$ and $(14)(25)(3)$ acting on the tree shown in Figure \ref{fig:perm_example}.

Let $S_n$ denote the permutation group of $[n] = \{1,2,\ldots,n\}$.
In this case, the group action just defined extends uniquely, linearly, to an action of probabilities on $S_n$ upon probabilities on (fat or thin) cladograms.
The action of the element $\frac{1}{n!}\sum_{\sigma \in S_n} \sigma$ will be of interest later on.  This element has the effect of applying a uniform random permutation to the leaf labels of a cladogram.

\subsection{Useful constructions on trees}
This section covers the {\em root join} operation on trees, and the set of splits of a tree.  The root join is used in the next chapter to define the alpha models.  The splits of a tree are used to calculate the probability of a given tree under these models.

Informal definitions are given first, followed by more rigorous definitions and proofs.

\begin{figure}[htbp]
    \begin{center}
    \resizebox{10cm}{!}{\includegraphics{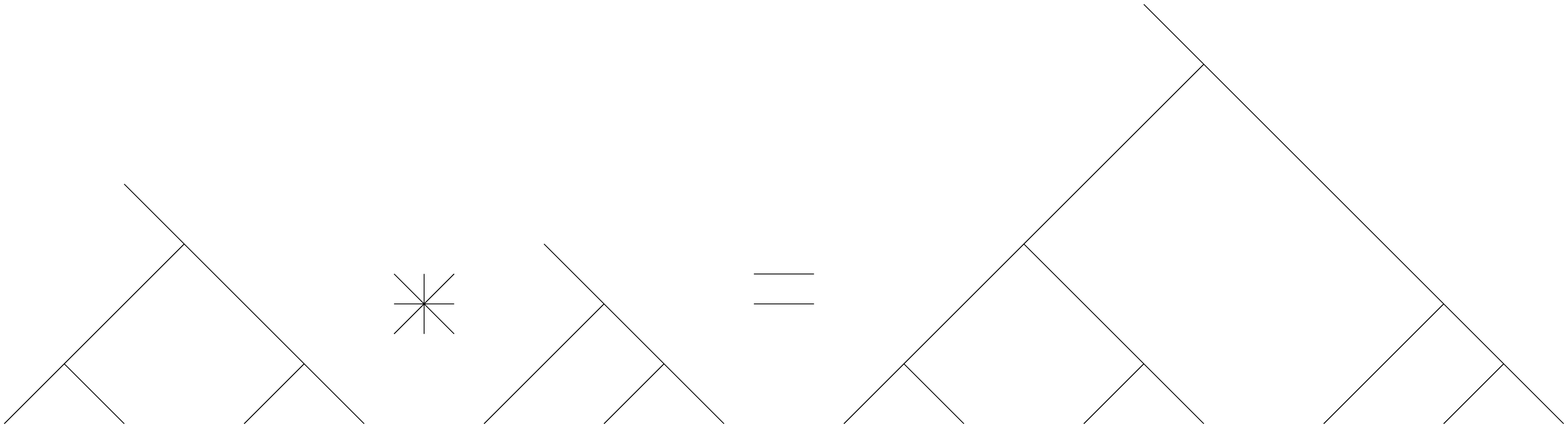}}\\
    \caption{\label{fig:root-join-example}Joining two trees at the root}
    \end{center}
\end{figure}

If $t$ is a fat rooted binary tree which has left subtree $t_1$ and right subtree $t_2$ then $t$ may be thought of as the tree formed by joining together $t_1$ and $t_2$ at their roots.  This is denoted $t_1*t_2=t$.  See Figure \ref{fig:root-join-example} for an example.

This construction also makes sense for thin (non-fat) trees and labeled trees, and this {\em root join} operation is preserved by the maps which forget orientation or leaf labels.  Every binary tree is the root join of two subtrees in this way.

For thin trees $t_1*t_2 = t_2*t_1$, but this is not true in general for fat trees.

\begin{figure}[htbp]
    \begin{center}
    \resizebox{10cm}{!}{\includegraphics{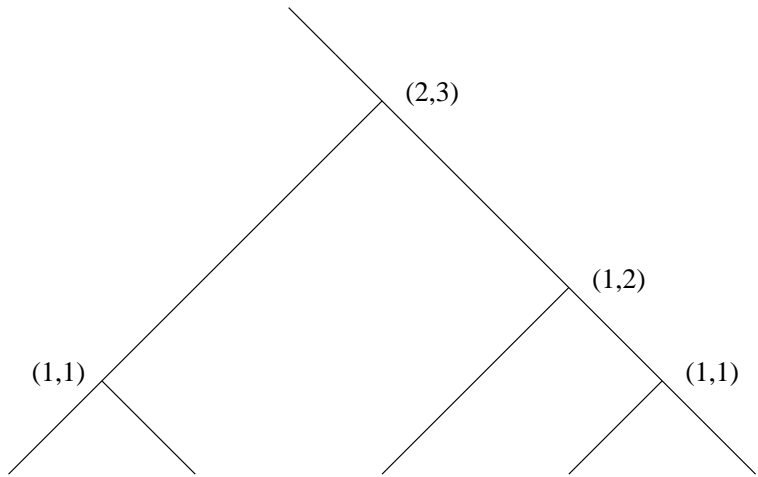}}\\
    \caption{\label{fig:splits-example}The splits of a tree}
    \end{center}
\end{figure}

{\em Splits} are defined as follows. If $t = t_1*t_2$ and $t_i$ has $n_i$ leaves then say that the first split of $t$ is the ordered pair $(n_1,n_2)$ (or the unordered pair if $t$ is a thin tree).  Similarly, each internal node has an associated split, as it is a branching point with some number of leaves below and to the left or right.  The multiset (set with multiplicity) of splits of a tree is useful for calculating the probability of a tree under certain classes of self-similar probabilities.  See Figure \ref{fig:splits-example} for an example.

\subsection{The subtree below an edge}
The definition of the subtree below an edge is useful in defining the root join operation.

\begin{figure}[htbp]
    \begin{center}
    \resizebox{10cm}{!}{\includegraphics{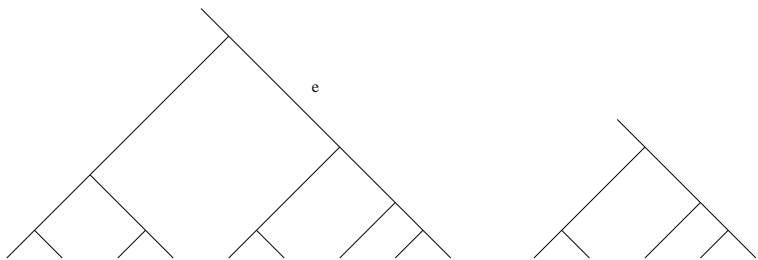}}\\
    \caption{\label{fig:subtree-below-an-edge}The subtree below an edge}
    \end{center}
\end{figure}

\begin{definition}
Given an edge $e$ of a rooted tree $t$, the {\em subtree of $t$ below edge $e$}, call it $s$, is a the rooted tree with 

\noindent vertex set $V$ comprising all vertices of $t$ which are descendents of both ends of $e$, 

\noindent edge set comprising all edges of $t$ which have both ends in $V$ and 

\noindent root vertex the end of $e$ which is closest to the root (and so of degree $1$ in the new tree).

\noindent Furthermore:

If $t$ is a partially labeled tree then so is $s$, with labeling function the restriction of the original labeling function to $V$.  Thus every vertex of $s$ is labeled exactly as it was in $t$.

If $t$ is a fat tree then so is $s$ and the orientation of every vertex of $s$ is the same as the orientation of that vertex in $t$, with the exception of the root of $s$ which has orientation the length-one cycle consisting of it's unique neighbor.  This is the only possible choice of orientation at the root vertex.
\end{definition}
See Figure \ref{fig:subtree-below-an-edge} for and example.

\begin{proposition}
The graph called the {\em subtree of $t$ below edge $e$} is indeed a rooted tree
\end{proposition}
\begin{proof}
First, show that the subtree of $t$ below edge $e$, call it $s$, is a tree.  Let $e = \{v_1,v_2\}$, with $v_1$ closer to the root of $t$, so that the root of $s$ is defined to be $v_1$.
Any path in $s$ is a path in $t$ and so there is at most one path between any two vertices of $s$.  On the other hand, every vertex of $s$ is either $v_1$ or a descendent of $v_2$ (ancestors and descendents referring to tree $t$).  Every non self-intersecting path from a descendent of $v_2$ to $v_2$ passes through descendents of $v_2$ only.  Since all descendents of $v_2$ lie in $s$ it follows that $s$ is connected.  Thus $s$ is a tree.  Finally, $v_1$ has degree $1$ in $s$, so $s$ is a rooted tree.
\end{proof}

\subsection{Joining two trees at the root}
This section contains the definition of the operation $\cdot*\cdot$ of joining two trees at the root.  This operation is used extensively in the definitions to come so several of its properties are examined in detail.

\begin{definition}
Given rooted trees $t_1$ and $t_2$, let the {\em root join of $t_1$ and $t_2$}, denoted $t_1 * t_2$, be the tree defined as follows:
\begin{itemize}
\item If $t_1$ is an empty tree then $t_1*t_2 = t_2$.  If $t_2$ is an empty tree then $t_1*t_2 = t_1$.
\item Otherwise, the tree $t_1*t_2$ includes vertices $r,v_0,v_1,v_2$ and edges $\{r,v_0\}$, $\{v_0,v_1\}$, $\{v_0,v_2\}$, such that $r$ is the root vertex and, for each $i$ in $\{1,2\}$, the subtree of $t_1*t_2$ below edge $\{v_0,v_i\}$ is isomorphic, via $f_i$, to $t_i$.
\end{itemize}

\noindent Furthermore:

If $t_1$ and $t_2$ are leaf-labeled trees then so is $t_1*t_2$, and the maps $f_1,f_2$ are isomorphisms of partially labeled trees.

If $t_1$ and $t_2$ are fat (oriented) trees then $t_1*t_2$ is a fat tree, the maps $f_1,f_2$ are isomorphisms of fat trees, the orientation at $r$ is the cycle $(v_0)$ and the orientation at $v_0$ is the cycle $(r,v_1,v_2)$.
\end{definition}

\begin{proposition}
Given rooted trees $t_1$ and $t_2$ as in the previous definition, the tree denoted $t_1*t_2$ exists and is uniquely defined up to isomorphism.
\end{proposition}
\begin{proof}
If $t_1$ or $t_2$ is the empty tree then $t_1*t_2$ is equal to either $t_2$ or $t_1$ and so exists and is uniquely defined.
Suppose now that $t_1$ and $t_2$ are non-empty trees.

Assume for the moment that $t_1$ and $t_2$ are thin, unlabeled rooted trees.

Let $t_i$ have vertex set $V_i$ and edge set $E_i$.  Without loss of generality, suppose that the vertex sets of $t_1$ and $t_2$ intersect at a single element, $v_0$, which is the root for both trees.  Let $r$ be an element not contained in $V_1$ or $V_2$.  This will represent the root of the new tree.

Let $V = V_1 \cup V_2 \cup \{r\}$ and $E = E_1 \cup E_2 \cup \{ \{v_0,r\} \}$.  Let $t$ be the graph with vertex set $V$ and edge set $E$.  Now show that $t$ has the properties required of $t_1*t_2$.

First, show that $t$ is a tree.  The graph $t$ is connected as there is a path from every vertex to the vertex $v_0$.  Now to show that there is a unique non self-intersecting path between any two vertices.  Note that $t_1$, $t_2$, and the graph $t_3$ with vertex set $V_3 = \{v_o,r\}$ and edge set $E_3 = \{\{v_0,r\}\}$ are all trees.

Any edge from a vertex in $V_i$ to a vertex in $V_j$, with $i\ne j$ must contain vertex $v_0$.
Thus, any path from a vertex in $V_i$ to a vertex in $V_j$, with $i\ne j$ must pass through $v_0$.
Thus, any non self-intersecting path from a vertex, $x$, in $V_i$ to a vertex, $y$, in $V_j$ must contain $v_0$ exactly once, with all vertices in the path before $v_0$ lying in $V_i$ and all those after $v_0$ lying in $V_j$.  The sub-path from $x$ to $v_0$ in non self-intersecting and lies entirely in $V_i$ and so is unique, since $(V_i,E_i)$ is a tree.  Similarly with the sub-path from $v_0$ to $y$.  Thus the non self-intersecting path from $x$ to $y$ must be unique.

Any non self-intersecting path from a vertex $x$ to $y$, both in $V_i$, must lie entirely in $V_i$.  Otherwise, if $v$ is any vertex in the path not lying in $V_i$ (and so not equal to $v_0$) then the sub-path from $x$ to $v$ passes through $v_0$ as does the sub-path from $v$ to $y$.  Thus $v_0$ appears twice on a non self-intersecting path, which is a contradiction.  Therefore, since the non self-intersecting path from $x$ to $y$ lies entirely in $V_i$ is must be unique, since $(V_i,E_i)$ is a tree.

Now show that $t$ has the required properties.  First, It contains the required vertices, $r$,$v_0$,$v_1$,$v_2$, and edges, $\{r,v_0,\}$,$\{v_0,v_1\}$,$\{v_0,v_2\}$, which are explicitly stated.  Second, by the construction of $t$, for each $i=1,2$, the subtree of $t$ below edge $\{v_0,v_i\}$ has vertex set $V_i$ and edge set $E_i$ and therefore is isomorphic to $t_i$.

Furthermore if $t_1$ and $t_2$ are fat rooted trees, with orientation functions $o_1$ an $o_2$, then let $t$ have orientation function $o$ defined by $o(r) = (v_0)$, $o(v_0) = (r,v_0,v_1)$ and $o(v) = o_i(v)$ for $v\in V_i \\ \{v_0\}$.  Thus $t$ satisfies the additional requirements on the orientation of $t_1*t_2$.

Note that vertex $v$ is a (non-root) leaf of $t$ if and only if  $v$ is a leaf of either $t_1$ or $t_2$.

Furthermore, if $t_1$ and $t_2$ are leaf-labeled trees then let $t$ be a leaf-labeled tree such that the leaf $v \in V_i \subset V$ of $t$ has the same label both as a vertex of $t$ and of $t_i$.  All other vertices of $t$ are unlabeled.  Thus $t$ satisfies the additional requirements on the labeling of $t_1*t_2$.

Next show that any two trees satisfying the definition of $t_1*t_2$ must be isomorphic.
Again, begin by assuming simply that $t_1$ and $t_2$ are thin unlabeled rooted trees.

Let $s_1$ and $s_2$ be trees which satisfy the requirements of $t_1*t_2$.
Therefore, $s_i$ contains vertices $r_i,v_{i0},v_{i1},v_{i2}$ and edges $\{r_i,v_{i0}\}$, $\{v_{i0},v_{i1}\}$, $\{v_{i0},v_{i2}\}$, such that $r_i$ is the root vertex of $s_i$ and, for each $j$ in $\{1,2\}$, the subtree of $t_1*t_2$ below edge $\{v_{i0},v_{ij}\}$ is isomorphic, via $f_{ij}$, to $t_j$.

Let $f$ be a map from the vertex set of $s_1$ to the vertex set of $s_2$ defined such that $f(r_1) = r_2$,
$f(v_{1j}) = v_{2j}$ for $j=0,1,2$ and if $v$ is a descendent of $v_{1j}$ then $f(v) = f_{2j}^{-1}f_{1j}(v)$.  This map is a bijection on vertices, sends the root to the root, and maps edges to edges, as does its inverse.  Thus is it an isomorphism of thin rooted trees.

If $t_1$ and $t_2$ are both fat trees, or both leaf-labeled trees, then $f$ is also an isomorphism of, respectively, fat trees or leaf-labeled trees.

Thus there is exactly one tree, up to isomorphism, satisfying the requirements of $t_1*t_2$.
\end{proof}

Note that the sum of the number of leaves in two rooted trees is the same as the number of leaves in the root join of these two trees.  In other words $|t_1*t_2| = |t_1| + |t_2|$.

Also, note that forgetting orientations before or after joining two trees at the root has the same effect.  The same is true for forgetting leaf-labelings.
Since this result is used often, it deserves a proposition.
\begin{proposition}
\label{prop:star-forget-commutes}
If $t_1$ and $t_2$ are rooted binary trees which are both fat then $F_o(t_1*t_2) = F_o(t_1)*F_o(t_2)$, and if $t_1$ and $t_2$ are rooted binary trees which are both leaf-labeled then $F_l(t_1*t_2) = F_l(t_1)*F_l(t_2)$
\end{proposition}
\begin{proof}
This follows directly from the definitions for the binary operator $\cdot *\cdot$ and the forgetful maps $F_o$, which forgets orientations of fat trees, and $F_l$ which forgets labelings of labeled trees.
\end{proof}

The following result shows that there is only one way, up to isomorphism, to write a binary tree as the root join of two non-empty trees.
\begin{lemma}
\label{lemma:star-unique-fact}
If $t$ is a non-empty rooted binary tree then: either $t$ has one (non-root) leaf; or $t = t_1 * t_2 $ for
a unique pair of non-empty trees $\{t_1,t_2\}$, and if $t$ is a fat tree then there is a unique ordered pair $(t_1,t_2)$ such that $t = t_1*t_2$.
\end{lemma}
\begin{proof}
Suppose that $t$ is a non-empty rooted binary has more than $1$ leaf.
Therefore $t$ has a root, $r$, which has a unique neighbor, $v_0$.  This vertex, $v_0$, has degree $3$ and so has two distinct neighbors, $v_1$ and $v_2$, which are not the root $r$.  If $t$ is a fat tree then choose $v_1,v_2$ so that the orientation at $v_0$ is $(r,v_1,v_2)$.  For $i=1,2$, let $t_i$ be the subtree of $t$ below edge $\{v_0,v_i\}$.  Thus the tree $t_1*t_2$ is exactly the tree $t$ (provided that the root vertex of $t_1*t_2$ is chosen to be the same element as the root vertex of $t$).

Suppose that $t$ is also equal to $t_3*t_4$.  By the definition of the operation $\cdot * \cdot$, $t_3$ is isomorphic to the subtree of $t$ below edge $\{v_0,v_i\}$ for some $i \in \{1,2\}$ and $t_4$ is isomorphic to the subtree of $t$ below the other edge $\{v_0,v_j\}$, $j\in \{1,2\}$ such that $i\ne j$.  Thus $t_3$ is isomorphic to one of $t_1$ or $t_2$, and $t_4$ is isomorphic to the other.

Furthermore, if $t$ is a fat tree then $t= t_3*t_4$ implies that $t_3$ is isomorphic to the subtree of $t$ below edge $\{v_0,v_1\}$, which is $t_1$, and so the ordering of the two trees is also unique.
\end{proof}

\subsection{Splits}
Now for the formal definition of the {\em splits} of a tree.
First, if $t$ is a binary rooted tree then let $|t|$ denote the number of leaves of $t$, also called the {\em size} of $t$.

\begin{definition}
\label{def:first-split}
Suppose that $t = t_1 * t_2$ for non-trivial fat, respectively thin, rooted binary trees $t_1$ and $t_2$.  Say that $t$ has {\em first split} $(|t_1|,|t_2|)$, respectively $\{|t_1|,|t_2|\}$.
\end{definition}
Lemma \ref{lemma:star-unique-fact} ensures that the first split is well defined.

\begin{definition}
\label{def:family-of-splits}
Define the {\em family of splits} of a fat (respectively thin) rooted binary tree $t$ inductively as follows:

$\text{splits}(t)$ is a multiset (a set with multiplicities) such that
\begin{itemize}
\item If $t$ is a one-leaf tree then $\text{splits}(t) = \emptyset$
\item If $t =  t_1*t_2$, for non-trivial $t_1,t_2$, then 
\\ $\text{splits}(t) = \text{splits}(t_2)\dot\cup\text{splits}(t_2)\dot\cup\{(|t_1|,|t_2|)\}$ for fat trees, and
\\ $\text{splits}(t) = \text{splits}(t_2)\dot\cup\text{splits}(t_2)\dot\cup\{\{|t_1|,|t_2|\}\}$ for thin trees.
\end{itemize}
\end{definition}
Again, Lemma \ref{lemma:star-unique-fact} ensures that this is well defined.

\noindent An equivalent non-recursive definition is:
\begin{definition}
Given a rooted binary tree $t$, let $E_2$ be the set of edges which do not contain a leaf.  Define the multiset of splits of the tree $t$ to be the union over edges $e \in E_2$ of the first split of the subtree of $t$ below $e$.
\end{definition}

Equivalence of these definitions is not proven here.

\section{The alpha models}
Now that the requisite constructions and definitions are at hand, the alpha models may be defined.

The alpha models are four parameterized sequences of probability measures, one for each of the four types of trees focused on here: tree shapes, fat tree shapes, cladograms and fat cladograms.  For each type of tree, the $n$-th element of the corresponding sequence is a probability measure on the set of trees of that type with exactly $n$ leaves.  Each alpha model has a single real parameter $\alpha \in [0,1]$.

Each of the four sequences is constructed in a similar manner to the others, using successive {\em alpha insertions} to build up each probability measure.  The four are related through the maps which forget orientation and leaf-labels.  Each also has two interesting properties, called Markovian self-similarity and deletion stability (also called sampling consistency).  These two properties are briefly mentioned below and properly defined in the following sections.

The alpha model on cladograms is perhaps of most practical interest.  It is also representative of all four models, and is now described.

{\em Alpha insertion} of a leaf labeled $k$ into a cladogram is performed as follows.  Give each leaf edge weight $1-\alpha$ and all other edges weight $\alpha$.  Choose an edge at random according to these weights and attach a new leaf edge to the middle of this edge.  Label the newly created leaf $k$.
\begin{figure}[htbp]
    \begin{center}
\ignore{    \resizebox{10cm}{!}{\includegraphics{alpha_weights.eps}}\\}
    \resizebox{10cm}{!}{\includegraphics{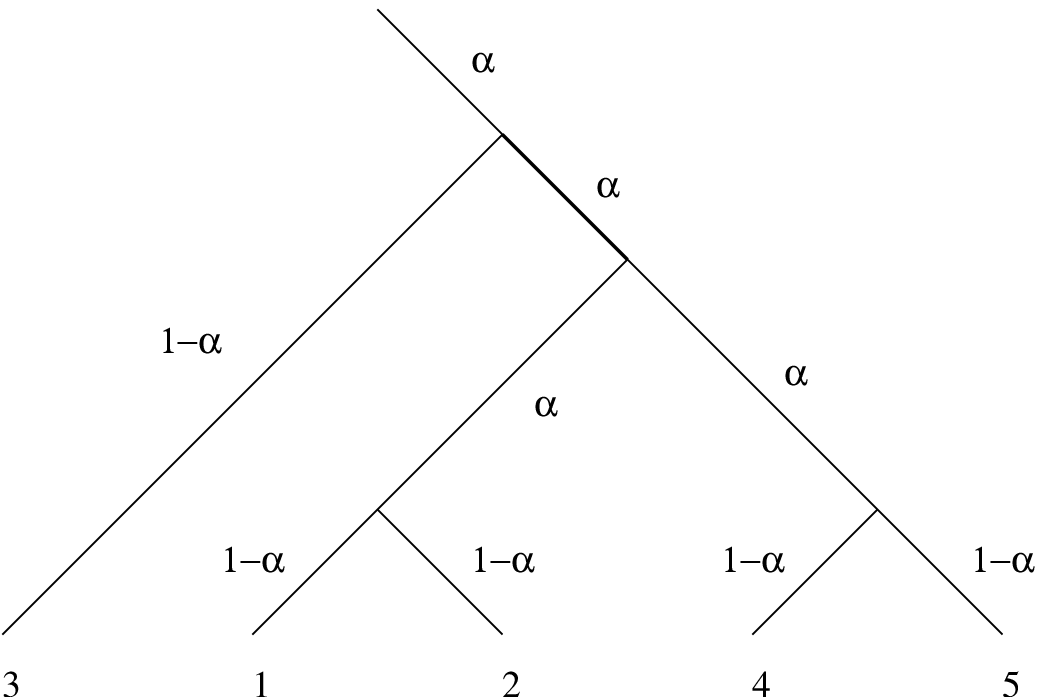}}\\
    \caption[Alpha insertion by edge weight]{\label{fig:weight-example2}The weight of a leaf edge is $1-\alpha$, the weight of an internal edge is $\alpha$}
    \end{center}
\end{figure}

\begin{figure}[htbp]
    \begin{center}
    \resizebox{10cm}{!}{\includegraphics{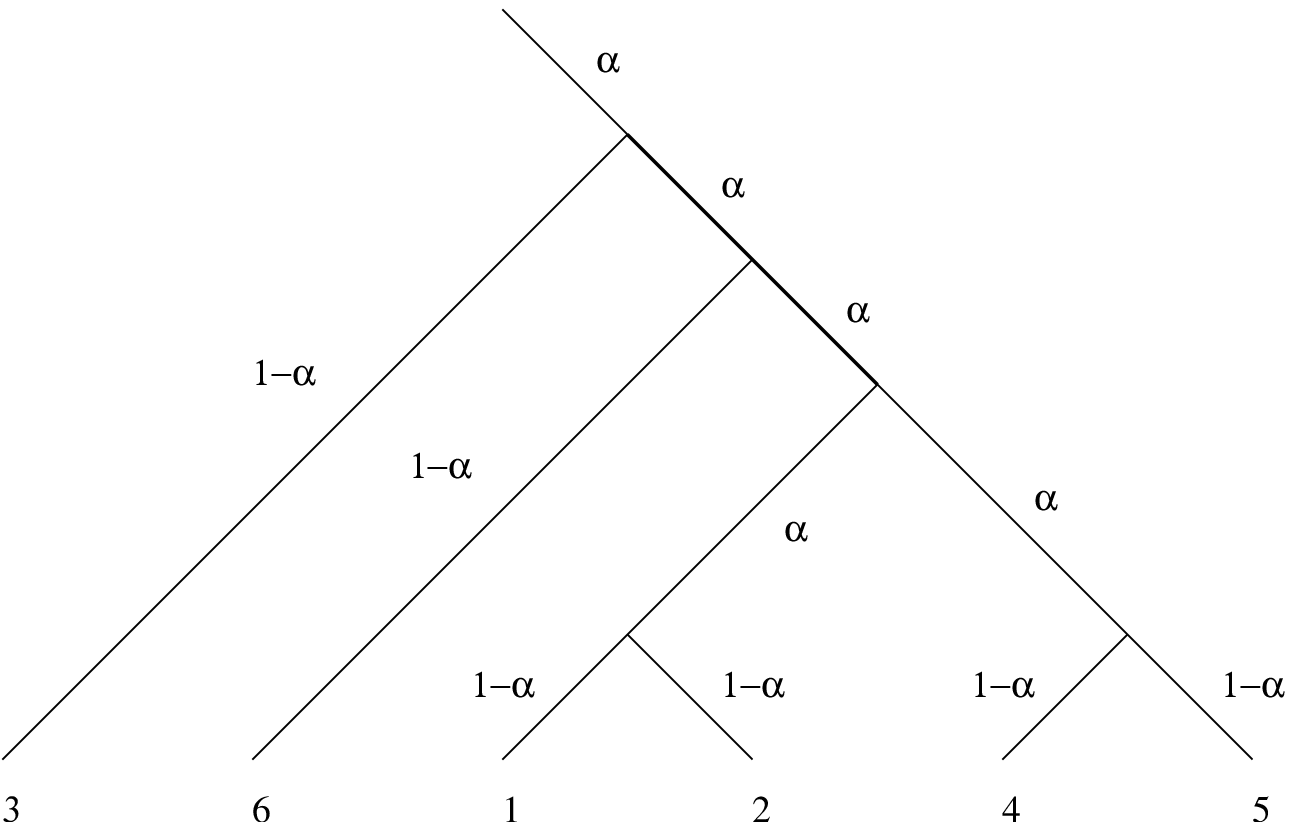}}\\
    \caption[Alpha insertion by edge weight]{The resulting cladogram with weights after inserting into the highlighted edge in Figure \ref{fig:weight-example2}}
    \end{center}
\end{figure}

A random cladogram with $n$ leaves from this model may be constructed as follows.
Take a rooted tree with a single leaf and label this leaf $1$.
Successively insert leaves labeled $2,3,\ldots,n$ into the tree according to the alpha insertion rule.
Once all of the leaves have been inserted, apply a uniform random permutation to the leaf labels.
Thus, the resulting distribution on cladograms is symmetric under permutation of leaf labels.

It is also {\em deletion stable} in the sense that if a random leaf is deleted (without loss of generality, the leaf with the largest label), then the resulting smaller tree is also distributed according to the alpha model with the same value of $\alpha$.  This implies {\em sampling consistency}: Given a random cladogram from the alpha distribution on cladograms with $n>k$ leaves, the shape of the subtree spanned by leaves $1,2,\ldots,k$ is distributed as the alpha distribution on cladograms with $k$ leaves.  In the case of unlabeled trees, a subtree spanned by $k$ randomly chosen leaves is distributed as an unlabeled alpha tree with $k$ leaves.
Deletion stability is described in more detail in Section \ref{section:deletion-stability}.

Another nice property of the alpha model is that it is {\em Markovian self-similar} (also called {\em Markov branching}).  This means that if the subtree below any edge has $k$ leaves, then the shape of this subtree is distributed as the alpha model on $k$ leaf trees, and is independent from the shape of the rest of the tree (conditional on there being exactly $k$ leaves below the given edge).  This is covered in Section \ref{section:mss}, in particular Proposition \ref{lemma:alpha-mss}.

Note that neither Markovian self-similarity nor sampling consistency (deletion stability) implies the other.

The alpha models are also the stationary distributions of certain Markov chains (one for each model).  These Markov chains 'project' onto each other via the forgetful maps.  These will be examined in later work.
\ignore{They are described in Section \ref{chapter:DAI-markov-chains} and their eigenvalues and mixing times examined in detail.}

In the formal definitions and proofs to follow, a recursive definition of alpha insertion is used.
Notice that if there are $n=n_1+n_2$ leaves total and $n_1$ leaves below one side of the first branch point then the probability that the new leaf is inserted in some edge down that branch is $\frac{n_1 - \alpha}{n-\alpha}$, the probability that it is inserted in some edge down the other side is $\frac{n_2 - \alpha}{n-\alpha}$ and the probability that it is inserted at the root edge is $\frac{\alpha}{n-\alpha}$.  This observation is the basis of the recursive definition of alpha insertion.

\subsection{Recursive definitions of alpha insertion}
The recursive definitions of alpha insertion for each type of tree are given below.

For the remainder of this section, let $s$ denote the one leaf binary rooted tree (tree shape), which is fat or thin as required by the context.
This tree has two vertices and a single edge from the root vertex to the non-root leaf.
Let $s_x$ be the one leaf binary rooted tree with leaf labeled $x$.

Let $|t|$ denote the number of leaves of a binary rooted tree $t$.

\begin{definition}
\label{def:alpha-insert-1}
Let $t$ be a fat binary rooted tree (fat tree shape).  Define $i_\alpha(t)$ as follows.
If $t$ has one leaf then define $i_\alpha(t) = \frac{1}{2} (t*s + s*t)$.  If not, then $t = t_1 * t_2$ for unique non-trivial $t_1$ and $t_2$.  In this case define
$$i_\alpha(t) = \frac{|t_1|-\alpha}{|t|-\alpha} i_\alpha(t_1) * t_2 + \frac{|t_2|-\alpha}{|t|-\alpha} t_1 * i_\alpha(t_2) + \frac{\alpha}{|t|-\alpha} \frac{1}{2}( s * t + t*s)$$
\end{definition}

\begin{definition}
If $t$ is a thin binary rooted tree then $i_\alpha(t)$ is given by exactly the same formulae as in the case of a fat tree.
\end{definition}
Uniqueness of the unordered pair $\{t_1,t_2\}$ and commutativity of the root join operation on thin trees ensures that $i_\alpha$ is well defined in this case.

\begin{definition}
\label{def:i_alpha_x}
Let $i_{\alpha,x}$ be defined for (fat or thin) leaf labeled rooted trees identically to $i_\alpha$ with the exception that the unlabeled single leaf tree $s$ is replaced everywhere with the labeled single leaf tree $s_x$.
\end{definition}

\begin{proposition}
\label{prop:insert-forget-commutes}
Alpha insertion commutes with forgetting orientation or leaf labels.
\end{proposition}
In other words, if $F_l$ is the function forgets leaf labels and $F_o$ is the function which forget orientations then
$$F_l( i_{\alpha,x}(t)) = i_\alpha(F_l(t))$$
$$F_o(i_{\alpha,x}(t)) = i_{\alpha,x}(F_o(t))$$
$$F_o(i_{\alpha}(t)) = i_{\alpha}(F_o(t))$$
for trees $t$ of the appropriate type.

\begin{proof}
That alpha insertion commutes with forgetting orientations, $F_o$, follows directly from the definitions.
The case of forgetting leaf labels follows by a simple induction.

For the initial case, $F_l(s_x) = s$ and the root join operation is preserved by the map which forgets labels (Proposition \ref{prop:star-forget-commutes}).  Thus, if $t$ is a single leaf tree then $F_l(i_{\alpha,x}(t)) = F_l( \frac{1}{2}(t*s_x + s_x * t) = \frac{1}{2}(F_l(t)*F_l(s_x) +  F_l(s_x)*F_l(t)) = \frac{1}{2}(F_l(t)*s +  s*F_l(t)) = i_\alpha(F_l(t)$.

For the inductive step, if $t$ is not a single-leaf tree then $t = t_1*t_2$ for non-trivial trees $t_1,t_2$.  Assume the statement is true for all trees smaller than $t$, with fewer leaves that is.
As $F_l$ respects $*$ it follows that
$F_l(i_{\alpha,x}(t))$ is equal to
$$
\frac{|t_1|-\alpha}{|t|-\alpha} F_l(i_{\alpha,x}(t_1)) * F_l(t_2) + \frac{|t_2|-\alpha}{|t|-\alpha} F_l(t_1) * F_l(i_{\alpha,x}(t_2)) + \frac{\alpha}{|t|-\alpha} \frac{1}{2}( F_l(s_x) * F_l(t) + F_l(t)*F_l(s_x))
$$
By the inductive assumption this is equal to
$$
\frac{|t_1|-\alpha}{|t|-\alpha} i_{\alpha}(F_l(t_1)) * F_l(t_2) + \frac{|t_2|-\alpha}{|t|-\alpha} F_l(t_1) * i_{\alpha}(F_l(t_2)) + \frac{\alpha}{|t|-\alpha} \frac{1}{2}( s * F_l(t) + F_l(t)*s))
$$
which is equal to $i_\alpha(F_l(t))$ as desired, since $F_l(t) = F_l(t_1*t_2) = F_l(t_1)*F_l(t_2)$.
\end{proof}

\subsection{Definitions of the alpha models}
The definitions of the alpha models for each of the four classes of trees are very similar.  Each involves successive alpha insertions, of labeled or unlabeled leaves, and then a final uniform randomization of leaf labels in the labeled cases.

The definitions of the alpha models depend upon a single variable, usually called alpha or $\alpha$, which lies in the range $[0,1]$.  Assume throughout that $\alpha$ is some fixed number.

Let $u_n = \frac{1}{n!}\sum_{\sigma \in S_n} \sigma$ be the uniform probability measure on the permutations of $[n] = \{1,2,\ldots,n\}$.

\begin{definition}
\label{def:alpha-model-ribbon-cladograms}\label{def:alpha-model-fat-cladograms}
The {\em alpha model on fat cladograms} is a sequence of probability measures $(P_n)_{n=1}^\infty$, such that $P_i$ is a probability measure on the set of fat cladograms with $n$ leaves, $P_0$ is the unique measure on the single fat cladogram with zero leaves (the empty tree), and for all integers $n\ge 1$
$$P_{n} = u_n i_{\alpha,n} \cdots i_{\alpha,2} i_{\alpha,1} P_0$$
\end{definition}
In other words, since $i_{\alpha,1} P_0=P_1$ is the unique measure on the single leaf tree, this definition says: start with the single leaf tree with leaf labeled $1$, alpha insert leaves labeled $2$ up to $n$ and then randomly permute the leaf labels.

\begin{proposition}
\label{prop:fat-alpha-model-recursion}
If $(P_n)_{n=1}^\infty$ is the alpha model on fat cladograms then for all integers $n\ge 1$
$$P_n = u_n i_{\alpha,n} P_{n-1}$$
\end{proposition}
\begin{proof}
The right hand side of the equation is equal to $u_n i_{\alpha,n} u_{n-1} i_{\alpha,n-1} \cdots i_{\alpha,1} P_0$.  Since alpha insertion does not depend on the position of the labels of $t$, it follows that $i_{\alpha,n} u_{n-1} t = \sigma i_{\alpha,n}t$, where $\sigma$ is the image of $u_{n-1}$ under the usual injection of $S_{n-1}$ into $S_n$.
Thus the right hand side is equal to $u_n \sigma i_{\alpha,n} i_{\alpha,n-1} \cdots i_{\alpha,1} P_0$ which is equal to $P_n$ as $u_n \sigma = u_n$ for any permutation $\sigma$.
\end{proof}

Define the {\em alpha model} on thin cladograms, fat tree shapes, and thin tree shapes to be the image of the alpha model on fat cladograms under the appropriate forgetful maps, $F_o$ and $F_l$ which forget orientations and leaf-labels respectively.  Specifically:

\begin{definition}
\label{def:alpha-model-general}
If $(P_i)_{i=1}^\infty$ is the alpha model on fat cladograms, $F_l$ is the function which forgets leaf labels and $F_o$ is the function which forgets orientations then define:

$(F_o(P_i))_{i=1}^\infty$ to be the alpha model on cladograms.

$(F_l(P_i))_{i=1}^\infty$ to be the alpha model on rooted binary fat trees.

$(F_o F_l(P_i))_{i=1}^\infty = (F_l F_o(P_i))_{i=1}^\infty$ to be the alpha model on rooted binary trees.
\end{definition}

\begin{proposition}
If $(P_i)_{i=1}^\infty$ is the alpha model on cladograms then
$$P_{n} = u_n i_{\alpha,n} P_{n-1} = u_n i_{\alpha,n} \cdots i_{\alpha,1} P_0$$
\end{proposition}
\begin{proposition}
If $(P_i)_{i=1}^\infty$ is the alpha model on (fat or thin) tree shapes then:
$$P_{n} = i_{\alpha} P_{n-1} = \underbrace{i_\alpha \cdots i_{\alpha}}_{n} P_0$$
\end{proposition}
\begin{proof}
Both of these propositions follow immediately from the previous two definition, and the fact that alpha insertion and the forgetful maps $F_o$ and $F_l$ 'commute' (Proposition \ref{prop:insert-forget-commutes}), and that $F_l(\sigma t) = F_l(t)$ for any fat or thin cladogram $t$ and any permutation, $\sigma$, of leaf labels.
\end{proof}

\subsection{Markovian self-similarity}
\label{section:mss}
Markovian self-similarity basically means that the subtree below any edge is picked from the distribution on trees of the correct size, independently of the rest of the tree.  It is also called {\em Markov branching} by Aldous in \cite{Aldous-1996}, as each branching happens independently of those above or on other paths from the root.

\begin{definition}
\label{def:MSS-1}
Let $(P_n)_{n=1}^{\infty}$ be a sequence of probability measures where $P_i$ is a probability on (fat) rooted binary trees with $n$ leaves.  Say that $(P_n)_{i=1}^\infty$ is {\em Markovian self-similar} if there exist real numbers $q(a,b)\ge 0$, for all integers $a,b \ge 1$, such that, for all integers $n \ge 2$, $\sum_{m=1}^{n-1} q(m,n-m) = 1$ and 
$$P_{n} = \sum_{m=1}^{n-1} q(m,n-m) P_{m}*P_{n-m}$$
\end{definition}

In other words, the trees below each child of the first branch-point are distributed independently from the same sequence of probabilities, conditional on the number of leaves they each have.

Call $q(\cdot,\cdot)$ the {\em conditional split distribution} of $(P_n)_{n=1}^{\infty}$.

Must $q$ be unique?
\begin{proposition}
\label{prop:unique-symmetric-q}
Suppose that such a $q$ exists, then in the case of fat trees $q$ is unique, and in the case of thin trees there is a unique symmetric $q$.
\end{proposition}
\begin{proof}
In the case of fat rooted binary trees, $t_1*t_2 = t_3*t_4$ if and only if $t_1 = t_3$ and $t_2 = t_4$.  In particular, the number of leaves of $t_1$ is equal to the number of leaves of $t_3$, so $\sum_{m=1}^{n-1} q_1(m,n-m) P_{m}*P_{n-m} = \sum_{m=1}^{n-1} q_2(m,n-m) P_{m}*P_{n-m}$ if and only if $q_1 = q_2$.

In the case of thin rooted binary trees, $t_1*t_2 = t_3*t_4$ if and only if $t_1 = t_3$ and $t_2 = t_4$, or $t_1 = t_4$ and $t_2 = t_3$.  In particular, the number of leaves of $t_1$ is equal to the number of leaves of either $t_3$ or $t_4$, so $\sum_{m=1}^{n-1} q_1(m,n-m) P_{m}*P_{n-m} = \sum_{m=1}^{n-1} q_2(m,n-m) P_{m}*P_{n-m}$ if and only if $q_1(m,n-m) + q_1(n-m,m) = q_2(m,n-m) + q_2(n-m,m)$.

Thus in the case of thin rooted binary trees there is a unique $q$ such that $q(a,b) = q(b,a)$.
\end{proof}

Conditional split distributions for thin trees are henceforth assumed to be symmetric in this way unless otherwise stated.

\begin{definition}
\label{def:MSS-2}
Similarly, a sequence $(P_n)_{n=1}^{\infty}$ of probability measures on (fat/thin) cladograms is called {\em Markovian self-similar} when the corresponding sequence on unlabeled trees, $(F_l(P_n))_{n=1}^{\infty}$, is Markovian self-similar.
\end{definition}

\subsection{Markovian self-similarity of the alpha models}

Define
$\Gamma_\alpha(n) = (n-1-\alpha)(n-2-\alpha)\cdots (2-\alpha)(1-\alpha)$
with $\Gamma_\alpha(1) = 1$.  Thus $\Gamma_0$ is the usual gamma function on the integers.

\begin{lemma}
\label{lemma:alpha-mss}
The four alpha models are all Markovian self-similar with the same conditional split distribution:
\begin{equation}
\label{eqn:alpha-split-distribution}
q_\alpha(a,b) = \frac{\Gamma_\alpha(a)\Gamma_\alpha(b)}{\Gamma_\alpha(a+b)}
\left( \frac{\alpha}{2}\binom{a+b}{a} + (1-2\alpha)\binom{a+b-2}{a-1} \right)
\end{equation}
\end{lemma}
\begin{proof}
By Definition \ref{def:alpha-model-general} and Proposition \ref{prop:star-forget-commutes} it suffices to prove that the alpha model on fat tree shapes is Markovian self-similar with the specified split distribution.

First, use induction to show that the first split of the alpha model on fat trees is distributed according to $q_\alpha$.

Recall that $P_{n+1} = i_\alpha P_n$ (Proposition \ref{prop:fat-alpha-model-recursion}), and that if 
$t_1$ has $a$ leaves and $t_2$ has $b$ leaves then tree $t_1 * t_2$ has first split $(a,b)$ (Definition \ref{def:first-split}).

By the formula for alpha insertion, $i_\alpha$, (Definition \ref{def:alpha-insert-1}), if $t$ is a fat tree shape with first split $(a,b)$ then $i_\alpha t$ is a fat tree with first split:
\begin{itemize}
\item$(a+1,b)$ with probability $\frac{a-\alpha}{n-\alpha}$
\item$(a,b+1)$ with probability $\frac{b-\alpha}{n-\alpha}$
\item$(1,a+b)$ with probability $\frac{\alpha}{(n-\alpha)}\frac{1}{2}$
\item$(a+b,1)$ with probability $\frac{\alpha}{(n-\alpha)}\frac{1}{2}$
\end{itemize}

To start the induction, note that for $n=2$ there is only one fat tree, and $q_\alpha(1,1)=1$ as it should.

Next, suppose $t$ is a random fat tree shape with $n$ leaves.
Show that if the first split, $(a,n-a)$, of $t$ is distributed as $q_\alpha(a,n-a)$ then the first split of $i_\alpha t$ is distributed as $q_\alpha(a,n+1-a)$.

In other words, show that $q_\alpha$ satisfies the following equations:
\begin{eqnarray}
\label{eqn:alpha-split-recursions}
q_\alpha(1,1)	& = & 1 \text{   , and for all $a,b > 1$:}
\nonumber
\\ q_\alpha(a,b)	& = & q_\alpha(a-1,b) \frac{a-1-\alpha}{a+b-1-\alpha}  + q_\alpha(a,b-1) \frac{b-1-\alpha}{a+b-1-\alpha}
\nonumber
\\ q_\alpha(1,b)& = & \frac{\alpha/2}{b-\alpha} + q_\alpha(1,b-1)\frac{b-1-\alpha}{b-\alpha}
\\ q_\alpha(a,1)& = & \frac{\alpha/2}{a-\alpha} + q_\alpha(a-1,1)\frac{a-1-\alpha}{a-\alpha}
\nonumber
\end{eqnarray}
This computation is omitted.

Thus, by induction, the first split of the alpha model satisfies $q_\alpha$.

Next, to show that the alpha model is Markovian self-similar.  In other words, show that if $(P_n)_{i=0}^\infty$ is the alpha model on  fat tree shapes then
$$P_{n} = \sum_{m=1}^{n-1} q_\alpha(m,n-m) P_{m}*P_{n-m}$$
This equation is true for $n=2$.
Suppose it is true for some $n$, then $P_{n+1} = i_\alpha P_n = \sum_{m=1}^{n-1} q_\alpha(m,n-m) i_\alpha(P_{m}*P_{n-m})$

Recall that if $t=t_1*t_2$ then
$i_\alpha(t) = \frac{|t_1|-\alpha}{|t|-\alpha} i_\alpha(t_1) * t_2 + \frac{|t_2|-\alpha}{|t|-\alpha} t_1 * i_\alpha(t_2) + \frac{\alpha}{|t|-\alpha} \frac{1}{2}( s * t + t*s)$ (Definition \ref{def:alpha-insert-1})
It follows that $i_\alpha(P_{m}*P_{n-m})$ is equal to
$$\frac{m-\alpha}{n-\alpha} P_{m+1} *P_{n-m} + \frac{n-m-\alpha}{n-\alpha} P_m * P_{n-m+1} + \frac{\alpha}{n-\alpha} \frac{1}{2}( P_{1} * P_n + P_n*P_1)$$

Thus $P_{n+1}$ is a linear combination of terms of the form $P_{m}*P_{n+1-m}$ and so is equal to $\sum_{m=1}^{n} q(m,n+1-m) P_{m}*P_{n-m}$ for some $q$.
Since $q$ is the distribution of the first split, by the arguement above it must be equal to $q_\alpha$.
Thus, the inductive step holds and the proposition is proven.
\end{proof}

The proof of lemma \ref{lemma:alpha-mss}, whilst perfectly correct, gives no indication as to how the formula was derived in the first place.
One possible derivation of the formula is sketched in the discussion below.

\begin{discussion}
First, recall the recurrence relations for the conditional split distributions $q_\alpha$, shown in Equation set \ref{eqn:alpha-split-recursions}.

Next, use a network flow argument to find a closed form solution for the recurrence equations.

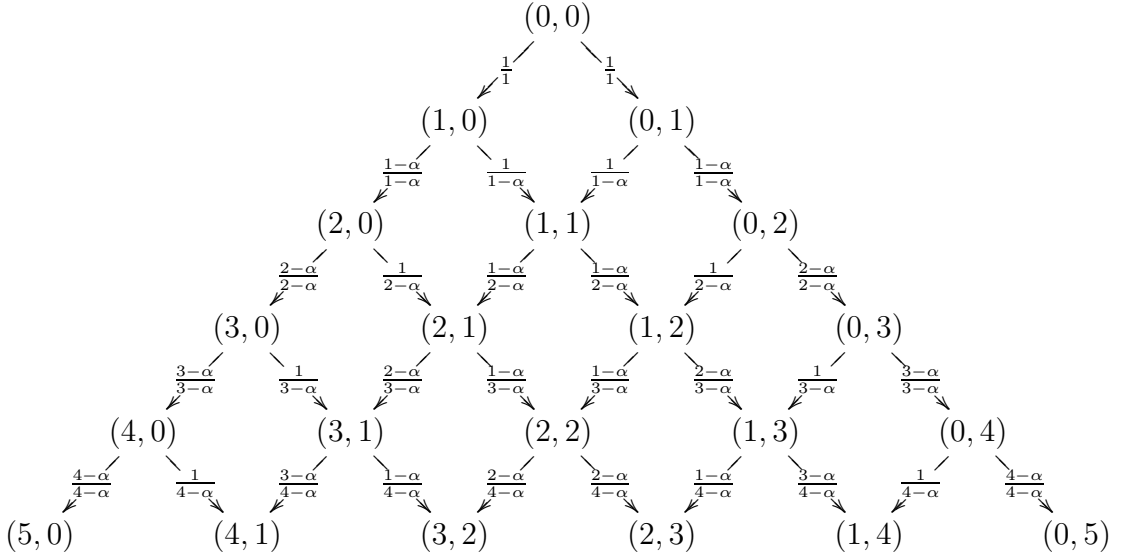
\begin{figure}[ht]
{\center
\xymatrix@+0.5pc@dr{
(0,0) \ar[d]|{\frac{1}{1}} \ar[r]|{\frac{1}{1}}& (0,1) \ar[d]|{\frac{1}{1-\alpha}} 
\ar[r]|{\frac{1-\alpha}{1-\alpha}}& (0,2) \ar[d]|{\frac{1}{2-\alpha}} \ar[r]|{\frac{2-\alpha}{2-\alpha}}& (0,3) 
\ar[d]|{\frac{1}{3-\alpha}} \ar[r]|{\frac{3-\alpha}{3-\alpha}}& (0,4) \ar[d]|{\frac{1}{4-\alpha}} 
\ar[r]|{\frac{4-\alpha}{4-\alpha}}& (0,5)\\
(1,0) \ar[d]|{\frac{1-\alpha}{1-\alpha}} \ar[r]|{\frac{1}{1-\alpha}}& (1,1) \ar[d]|{\frac{1-\alpha}{2-\alpha}} 
\ar[r]|{\frac{1-\alpha}{2-\alpha}}& (1,2) \ar[d]|{\frac{1-\alpha}{3-\alpha}} \ar[r]|{\frac{2-\alpha}{3-\alpha}}& 
(1,3) \ar[d]|{\frac{1-\alpha}{4-\alpha}} \ar[r]|{\frac{3-\alpha}{4-\alpha}}& (1,4)\\
(2,0) \ar[d]|{\frac{2-\alpha}{2-\alpha}} \ar[r]|{\frac{1}{2-\alpha}}& (2,1) \ar[d]|{\frac{2-\alpha}{3-\alpha}} 
\ar[r]|{\frac{1-\alpha}{3-\alpha}}& (2,2) \ar[d]|{\frac{2-\alpha}{4-\alpha}} \ar[r]|{\frac{2-\alpha}{4-\alpha}}& 
(2,3) \\
(3,0) \ar[d]|{\frac{3-\alpha}{3-\alpha}} \ar[r]|{\frac{1}{3-\alpha}}& (3,1) \ar[d]|{\frac{3-\alpha}{4-\alpha}} 
\ar[r]|{\frac{1-\alpha}{4-\alpha}}& (3,2) \\
(4,0) \ar[d]|{\frac{4-\alpha}{4-\alpha}} \ar[r]|{\frac{1}{4-\alpha}}& (4,1) \\
(5,0)
}
\caption[A flow network]{A flow network for the split distribution of the alpha model}
\label{fig:triangle-1}
}
\end{figure}

Think of the above triangular diagram (Figure \ref{fig:triangle-1}) as a directed flow network, with as yet unspecified sources and sinks.  The labels on each edge are the multiplying factor applied to the flow out of the starting vertex before it is added to the ending vertex.
Now choose the sources/sinks so that the net flow through each vertex $(a,b)$, for $a,b\ge 1$, is the conditional probability of the first split being $(a,b)$ when the tree has $a+b$ leaves total.

Notice that if this is true for one line $(1,n-1),\ldots,(n-1,1)$, then the contributions to the next line will be just as in Equations \ref{eqn:alpha-split-recursions} except for the $\frac{\alpha/2}{n-\alpha}$ contribution to $(1,n+1)$ and $(n+1,1)$.  This missing contribution should come from a flow of $\alpha/2$ through $(0,n)$ and $(n,0)$.

This implies that the node $(0,0)$ should be a source with inflow $\frac{\alpha}{2}$, so that the flows from $(0,n)$ to $(1,n)$ and $(n,0)$ to $(n,1)$ are $\frac{\alpha}{2}$ times $\frac{1}{n-\alpha}$ as needed.
This then implies that the node $(1,1)$ must be a source with enough inflow so that the total inflow is $1$, since $q_\alpha(1,1)=1$.  Thus it must be a source with inflow $\frac{1-2\alpha}{1-\alpha}$.  By the argument above, these are all the sources needed.
Figure \ref{fig:triangle-2} shows the network with the total flow into each node.

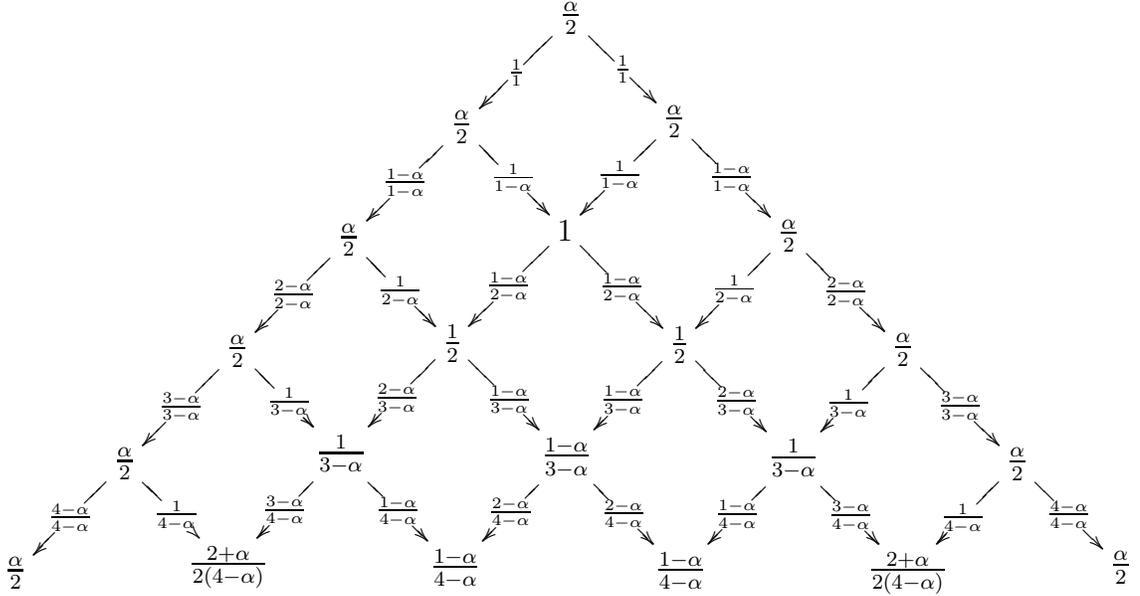
\begin{figure}[ht]
{\center
\xymatrix@+0.5pc@dr{
\frac{\alpha}{2} \ar[d]|{\frac{1}{1}} \ar[r]|{\frac{1}{1}}& \frac{\alpha}{2} \ar[d]|{\frac{1}{1-\alpha}} 
\ar[r]|{\frac{1-\alpha}{1-\alpha}}& \frac{\alpha}{2} \ar[d]|{\frac{1}{2-\alpha}} \ar[r]|{\frac{2-\alpha}{2-\alpha}}& \frac{\alpha}{2} 
\ar[d]|{\frac{1}{3-\alpha}} \ar[r]|{\frac{3-\alpha}{3-\alpha}}& \frac{\alpha}{2} \ar[d]|{\frac{1}{4-\alpha}} 
\ar[r]|{\frac{4-\alpha}{4-\alpha}}& \frac{\alpha}{2}\\
\frac{\alpha}{2} \ar[d]|{\frac{1-\alpha}{1-\alpha}} \ar[r]|{\frac{1}{1-\alpha}}& 1 \ar[d]|{\frac{1-\alpha}{2-\alpha}} 
\ar[r]|{\frac{1-\alpha}{2-\alpha}}& \frac{1}{2} \ar[d]|{\frac{1-\alpha}{3-\alpha}} \ar[r]|{\frac{2-\alpha}{3-\alpha}}& 
\frac{1}{3-\alpha} \ar[d]|{\frac{1-\alpha}{4-\alpha}} \ar[r]|{\frac{3-\alpha}{4-\alpha}}& \frac{2+\alpha}{2(4-\alpha)}\\
\frac{\alpha}{2} \ar[d]|{\frac{2-\alpha}{2-\alpha}} \ar[r]|{\frac{1}{2-\alpha}}& \frac{1}{2} \ar[d]|{\frac{2-\alpha}{3-\alpha}} 
\ar[r]|{\frac{1-\alpha}{3-\alpha}}& \frac{1-\alpha}{3-\alpha} \ar[d]|{\frac{2-\alpha}{4-\alpha}} \ar[r]|{\frac{2-\alpha}{4-\alpha}}& 
\frac{1-\alpha}{4-\alpha} \\
\frac{\alpha}{2} \ar[d]|{\frac{3-\alpha}{3-\alpha}} \ar[r]|{\frac{1}{3-\alpha}}& \frac{1}{3-\alpha} \ar[d]|{\frac{3-\alpha}{4-\alpha}} 
\ar[r]|{\frac{1-\alpha}{4-\alpha}}& \frac{1-\alpha}{4-\alpha} \\
\frac{\alpha}{2}  \ar[d]|{\frac{4-\alpha}{4-\alpha}} \ar[r]|{\frac{1}{4-\alpha}}& \frac{2+\alpha}{2(4-\alpha)} \\
\frac{\alpha}{2}
}
\caption[Filled in network]{The flow at each node of the network}
\label{fig:triangle-2}
}
\end{figure}

Finally, notice that all paths between any two nodes have the same product.
A path from $(0,0)$ to $(a,b)$ has weight 
$\frac{\Gamma_\alpha(a)\Gamma_\alpha(b)}{\Gamma_\alpha(a+b)}$, and one from $(1,1)$ has weight 
$\frac{\Gamma_\alpha(a)\Gamma_\alpha(b)}{\Gamma_\alpha(a+b)} (1-\alpha)$.  Also, note that there are $\binom{a+b}{a}$ possible paths from $(0,0)$ to $(a,b)$ and $\binom{a+b-2}{a-1}$ possible paths from $(1,1)$ to $(a,b)$.

Summing the inflow by path and source now gives the stated formula for $q_\alpha$.

\end{discussion}

\subsection{Calculating the probability of a tree}
\label{sec:probability-of-trees}
This section gives a simple method for calculating the probability of a tree under a Markovian self-similar sequence of probabilities, such as the alpha model.  Examples of the probabilities of small tree shapes are worked out.

First consider fat tree shapes (unlabeled rooted binary fat trees).

\begin{proposition}
\label{prop:tree-prob-fat-shapes}
Suppose that $(P_i)_i$ is a sequence of probabilities on fat tree shapes which is Markovian self-similar, with conditional split distributions given by $q$.  If $t$ is a tree with $n$ leaves whose family of splits is $F$ then
$$P_n(t) = \prod_{(a,b)\in F} q(a,b)$$
\end{proposition}
\begin{proof}
The statement is true for the single tree with one leaf, and for the single tree with two leaves.  For the pedantic, when $n=0$ the empty product is $1$ which is equal to the probability of the empty tree.

If the tree $t$ has at least $2$ leaves it may be written as $t = t_1*t_2$ and so $F = \text{splits}(t) = \{(a_1,b_1)\} \cup \text{splits}(t_1) \cup \text{splits}(t_2)$ (as a union of multisets).
The probability that random tree $t'$ has first split $(a_1,b_1)$ is $q(a_1,b_1)$.  Conditional on this, the probability that $t' = t= t_1*t_2$ is $P_a(t_1)P_b(t_2)$ (by the definition of Markovian self-similarity).  By induction this is $\prod_{(a,b) \in \text{splits}(t_1)} q(a,b) \prod_{(a,b) \in \text{splits}(t_1)} q(a,b)$.  Thus the probability of tree $t$ is $P_n(t) = \prod_{(a,b)\in F} q(a,b)$ as desired.
\end{proof}

If $q$ is a split distribution, then define $\hat q\{a,b\} = q(a,b) + q(b,a)$ if $a\ne b$ and $\hat q\{a,a\} = q(a,a)$.
\begin{proposition}
\label{prop:tree-prob-thin-shapes}
Suppose that $(P_i)_i$ is a sequence of probabilities on thin tree shapes which is Markovian self-similar, with conditional split distributions given by $q$.
If $t$ is an unlabeled thin rooted binary tree with $n$ leaves whose family of splits is $F$ then
$$P_n(t) = \prod_{(a,b)\in F} \hat q\{a,b\}$$
\end{proposition}
\begin{proof}
This proof is almost identical to that above.  In this case
the probability that $t$ has first split $\{a,b\}$ is $\hat q\{a,b\}$.
\end{proof}

For fat cladograms (labeled fat rooted trees):
\begin{corollary}
\label{prop:tree-prob-fat-clads}
Suppose that $(P_i)_i$ is a Markovian self-similar sequence of probabilities on fat cladograms such that if $F_l(t_1) = F_l(t_2)$ then $P_n(t_1) = P_n(t_2)$. (In other words, any two cladograms with the same shape have the same probability.)
Then if $t$ is a fat cladogram with $n$ leaves:
$$P_n(t) = \frac{1}{n!}\prod_{(a,b)\in \text{splits}(t)} q(a,b)$$
\end{corollary}
\begin{proof}
By Definition \ref{def:MSS-2} the family $(F_l(P_i))_i$ is a sequence of Markovian-self-similar probabilities on fat tree shapes (unlabeled fat rooted binary trees).
Since $F_l(t_1) = F_l(t_2)$ implies $P_n(t_1) = P_n(t_2)$ and the pre-image of any fat tree shape under $F_l$ has size $n!$ it follows that
$P_n(t) = \frac{1}{n!}F_l(P_n)(F_l(t))$.  Since the map $F_l$ does not change the family of splits of a tree it now follows that
$P_n(t) = \frac{1}{n!}\prod_{(a,b)\in \text{splits}(t)} q(a,b)$, as desired.
\end{proof}

\begin{figure}[htbp]
    \begin{center}
    \resizebox{10cm}{!}{\includegraphics{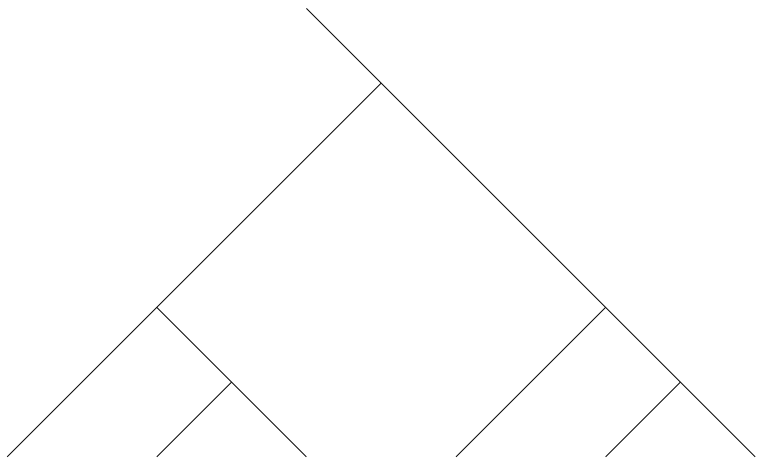}}\\
    \caption[Symmetric branch points]{\label{fig:symmetric-branchpoint-example}The first branch point of this tree is symmetric}
    \end{center}
\end{figure}

Before proceeding to the case of cladograms, a lemma is needed.
Say that a branch point is symmetric if the subtrees below each child edge are equal to each other.
See Figure \ref{fig:symmetric-branchpoint-example} for example.

\begin{lemma}
\label{lemma:number-of-clads-per-shape}
If $t$ is a tree shape with $n$ leaves and $k$ symmetric branch points then the number of cladograms with shape $t$ is equal to $n!/2^k$.
\end{lemma}
\begin{proof}
The symmetric group on $[n]=\{1,2,\ldots,n\}$ acts transitively on the set of cladograms with shape $t$.  The aim is now to show that the number of permutations which fix any cladogram with shape $t$ is $2^k$.  The lemma follows immediately from this.

Proceed by induction.
The statement to be proven is that a rooted binary tree with distinctly labeled leaves and $k$ symmetric branch points is fixed by $2^k$ permutations of its leaf labeling set.

First, this is trivially true for a labeled tree with $1$ leaf.

Suppose that the lemma is true for all trees smaller than $t$.

Suppose that the first branch-point of $t$ is not symmetrical (the easy case).
Then $t=t_1*t_2$ for distinct $t_1$ and $t_2$ and so any permutation which fixes $t$ must fix the set of leaves of $t_1$ and the set of leaves of $t_2$.  Thus the group fixing $t$ is the direct product of the group fixing $t_1$ and the group fixing $t_2$.  If $t_1$ has $k_1$ symmetric branch points and $t_2$ has $k_2$ symmetric branch points then $t$ has $k=k_1+k_2$ symmetric branch points.  Therefore, by the inductive assumption, the group fixing $t$ has size $2^k = 2^{k_1}2^{k_2}$.

Suppose that the first branch point of $t$ is symmetrical (the hard case).
Then $t = t_1*t_2$ where $t_1$ and $t_2$ have the same shape.  As they have isomorphic shapes, $t_1$ and $t_2$ both have the same number of symmetric branch points, say $k_1$.
Thus $t$ has $k = 1 + 2k_1$ symmetric branch points.

Now, every permutation which fixes $t$ must fix the unordered partitioning of leaf labels into those of one subtree and those of the other.
Thus, any permutation which fixes $t$ must either swap the two parts or not.  In each case, by the inductive assumption there are then $2^{k_1}$ distinct ways to permute the elements of each part without changing the cladogram.
Thus the order of the group fixing $t$ is $2 \times 2^{k_1} 2^{k_1} = 2^k$ as desired.
\end{proof}

Appendix \ref{chapter:tree-shapes-up-to-7} contains a list of all tree shapes with up to $7$ leaves, along with the number of cladograms of each shape.

So finally:
\begin{proposition}
\label{prop:tree-prob-thin-clads}
Suppose that $(P_i)_i$ is a Markovian self-similar sequence of probabilities on cladograms with split distribution $q$ such that if $F_l(t_1) = F_l(t_2)$ then $P_n(t_1) = P_n(t_2)$.  (In other words, any two cladograms with the same shape have the same probability.)
Then if $t$ is a cladogram with $n$ leaves and $k$ symmetric branch points:
$$P_n(t) = \frac{2^k}{n!}\prod_{(a,b)\in \text{splits}(t)} \hat q\{a,b\}$$
\end{proposition}
\begin{proof}
Similarly to the previous proof: by Definition \ref{def:MSS-2} the sequence $(F_l(P_i))_i$ is a sequence of Markovian self-similar probabilities on tree shapes (unlabeled rooted binary trees).
By Lemma \ref{lemma:number-of-clads-per-shape} the number of cladograms with the same shape as $t$ (ie such that  $F_l(t') = F_l(t)$) is $\frac{n!}{2^k}$ where $k$ is the number of equal splits of $t$.

Since $F_l(t_1) = F_l(t_2)$ implies $P_n(t_1) = P_n(t_2)$ and the pre-image of $F_l(t)$ under $F_l$ has size $\frac{n!}{2^k}$ it follows that
$P_n(t) = \frac{2^k}{n!}F_l(P_n)(F_l(t))$.  Since the map $F_l$ preserves the family of splits of a tree and the split distribution of a Markovian self-similar sequence, this gives:
$P_n(t) = \frac{2^k}{n!}\prod_{(a,b)\in \text{splits}(t)} q\{a,b\}$ as desired.
\end{proof}

\subsection{The probability of a tree under the alpha model}
\label{subsec:prob-of-trees-under-alpha}
There is little more to say in the special case of the alpha model.  Since the four alpha models satisfy the conditions, respectively, of Propositions \ref{prop:tree-prob-fat-shapes}, \ref{prop:tree-prob-thin-shapes}, \ref{prop:tree-prob-fat-clads} and \ref{prop:tree-prob-thin-clads}, the probability of a tree (of the appropriate type) under one of these models is given by the formulae in those propositions.

\begin{figure}[htbp]
    \begin{center}
    \resizebox{10cm}{!}{\includegraphics{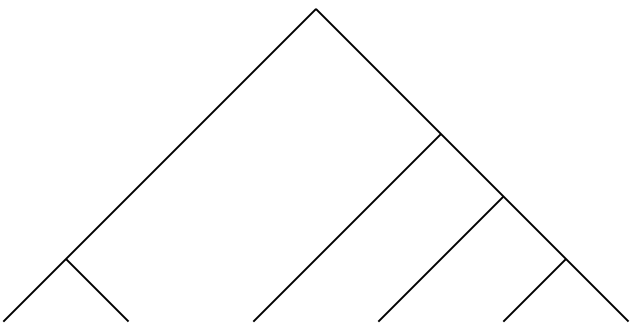}}\\
    \caption[A tree shape]{	\label{fig:example-treeshape-for-prob}
A tree shape with probability $\frac{2(1-\alpha)(8-\alpha)}{(5-\alpha)(4-\alpha)}$}
    \end{center}
\end{figure}

For example, the tree shape in Figure \ref{fig:example-treeshape-for-prob} has family of splits
$$\{ \{4,2\},\{1,1\},\{1,3\},\{2,1\},\{1,1\}\}$$

Therefore, by Proposition \ref{prop:tree-prob-thin-shapes}, the probability of this tree shape under the alpha model on tree shapes is 
$$ \prod_{(a,b)\in \{ \{4,2\},\{1,1\},\{1,3\},\{2,1\},\{1,1\}\} } \hat q\{a,b\}$$
Using Equation \ref{eqn:alpha-split-distribution} for the split distribution of the alpha model, and recalling that 
$\hat q\{a,b\} = q(a,b) + q(b,a)$ if $a\ne b$ and $\hat q\{a,a\} = q(a,a)$, this is equal to:
$$\frac{(1-\alpha)(8-\alpha)}{(5-\alpha)(4-\alpha)} \times 1 \times \frac{2}{3-\alpha} \times 1 \times 1$$
which simplifies to
$$\frac{2(1-\alpha)(8-\alpha)}{(5-\alpha)(4-\alpha)(3-\alpha)}$$

Appendix \ref{chapter:tree-shapes-up-to-7} contains a list of all tree shapes with up to $7$ leaves, along with the probability of each under the alpha model.

\subsection{Deletion stability}
\label{section:deletion-stability}
This section addresses the definition of {\em deletion stability}, and provides proofs that the alpha models have this property.

Informally, deletion stability on a sequence of probabilities $(P_i)_{i=0}^\infty$ on fat or thin cladograms means that picking a random cladogram with $n$ leaves from $P_n$ and deleting leaf $n$ gives a random cladogram with $n-1$ leaves distributed as $P_{n-1}$.
Similarly, a sequence of probabilities $(P_i)_{i=0}^\infty$ on fat or thin tree shapes is deletion stable if picking a random tree shape with $n$ leaves from $P_n$ and deleting a random leaf gives a random tree shape with $n-1$ leaves distributed as $P_{n-1}$.

The formal definition of deletion stability requires a formal definition of these deletions.

Let $D$ be the function which deletes a random leaf of a binary rooted tree.
A recursive definition of $D$ is given as this form is most convenient for the proofs which follow.  

For a tree shape or cladogram $t$, let $|t|$ denote the number of leaves of $t$, also called the size of $t$.
\begin{definition}
Let $t$ be a fat or thin tree shape.
\begin{itemize}
\item If $t$ is the empty tree then so is $D(t)$.
\item If $t$ has one leaf then $D(t)$ is the empty tree shape.
\item If $t$ has more than one leaf then $t=t_1*t_2$ for non-empty $t_1$,$t_2$, so let
$$D(t) = \frac{|t_1|}{|t_1|+|t_2|}D(t_1)*t_2 + \frac{|t_2|}{|t_1|+|t_2|}t_1*D(t_2)$$
\end{itemize}
\end{definition}

\begin{proposition}
The function $D$ is well defined.
\end{proposition}
\begin{proof}
If $t$ is a fat tree with more than one leaf then, by Lemma \ref{lemma:star-unique-fact}, $t=t_1*t_2$ for a unique pair of non-empty trees, $(t_1,t_2)$.  Thus $D(t)$ is well defined.  If $t$ is a thin tree then, by Lemma \ref{lemma:star-unique-fact}, $t=t_1*t_2=t_2*t_1$ for a unique set of two non-empty trees, $\{t_1,t_2\}$.  Since the root join operation is commutative for thin tree shapes 
$D(t_1*t_2) = \frac{|t_1|}{|t_1|+|t_2|}D(t_1)*t_2 + \frac{|t_2|}{|t_1|+|t_2|}t_1*D(t_2)
= \frac{|t_2|}{|t_1|+|t_2|}D(t_2)*t_1 + \frac{|t_1|}{|t_1|+|t_2|}t_2*D(t_1) = D(t_2*t_1)$ and so $D(t)$ is well defined.
\end{proof}

The next obvious result is that the operation $D$ is respected by the map, $F_o$, which forgets vertex orientations: taking fat tree shapes to thin tree shapes.
\begin{proposition}
\label{prop:delete-forget-orientation-commute}
If $t$ is a fat tree shape then $D(F_o(t)) = F_o(D(t))$.
\end{proposition}
\begin{proof}
If $t$ is a fat tree shape with one leaf then $F_o(t)$ is a thin tree shape with one leaf, $D(t)$ is the empty fat tree shape, and so both $F_o(D(t))$ and $D(F_o(t))$ are the empty thin tree shape.

Suppose that $t$ has more than one leaf and that the statement is true for all fat tree shapes with fewer leaves than $t$.
Now by Lemma \ref{lemma:star-unique-fact}, $t=t_1*t_2$ for a unique pair of non-empty trees $(t_1,t_2)$ and by Proposition \ref{prop:star-forget-commutes} $F_o(t) = F_o(t_1)*F_o(t_2)$.  Thus 
$$F_o(D(t)) = F_o\left( \frac{|t_1|}{|t_1|+|t_2|}D(t_1)*t_2 + \frac{|t_2|}{|t_1|+|t_2|}t_1*D(t_2) \right)$$
$$ = \frac{|t_1|}{|t_1|+|t_2|}F_o(D(t_1))*F_o(t_2) + \frac{|t_2|}{|t_1|+|t_2|}F_o(t_1)*F_o(D(t_2))$$
This is equal to $D(F_o(t))$, since by the inductive assumption $F_o(D(t_1)) = D(F_o(t_1))$ and $F_o(D(t_2)) = D(F_o(t_2))$, and the map $F_o$ leaves the number of leaves of a tree unchanged.  The result follows by induction.
\end{proof}

Now for the case of labeled trees.  The following is the definition of a function, $D_x$, which deletes every leaf labeled $x$.
\begin{definition}
\label{def:delete-x}
Let $t$ be a fat or thin labeled binary rooted tree.
\begin{itemize}
\item If $t$ is the empty tree then $D_x(t)$ is also the empty tree.
\item In the case where $t$ has one leaf: if this leaf is labeled $x$ then $D_x(t)$ is the empty tree, otherwise $D_x(t)=t$.
\item In the case where $t$ has more than one leaf: $t=t_1*t_2$ so define: $D_x(t) = D_x(t_1)*D_x(t_2)$.
\end{itemize}
\end{definition}
Again, Lemma \ref{lemma:star-unique-fact} guarantees that $D_x$ is well defined.

It is now shown that, for a random tree picked from a probability on cladograms invariant under permutation of leaf labels (such as the alpha model on cladograms), first deleting a specified leaf and then forgetting leaf lables is the same as first forgetting leaf labels and then deleting a random leaf.

\begin{proposition}
\label{prop:delete-forget-commute}
If $P$ is a probability on leaf-labeled (fat or thin) rooted binary trees with $n$ leaves which is invariant under any permutation of its leaf labels and such that all trees with positive probability have every leaf uniquely labeled and have a leaf labeled $x$, then $F_l(D_x(P)) = D(F_l(P))$.
\end{proposition}

\begin{proof}
First for the fat case.
If $n=1$ then a random tree $t$ from $P$ has one leaf, and this leaf is labeled $x$ and $D_x(t)$ is the empty tree.  Thus $F_l(D_1(P)) = D(F_l(P))$, the unique probability on the set of tree shapes with $0$ leaves, ie the empty tree.

Suppose that the result is true for all trees with less than $n$ leaves.

Let $t$ be a random tree picked from $P$, conditioned such that $t$ has first split $(a,b)$ (respectively $\{a,b\}$ in the thin case).  This implies that $t = t_1*t_2$ for a random pair of trees $(t_1,t_2)$ such that $|t_1|=a$ and $|t_2|=b$.
Now $D_x(t) = D_x(t_1)*D_x(t_2)$.   The leaf labeled $x$ is a leaf of $t_1$ with probability $\frac{|t_1|}{|t_1|+|t_2|}$, and in this case $x$ is not a leaf of $t_2$ and so $D_x(t) = D_x(t_1)*t_2$.  With probability $\frac{|t_1|}{|t_1|+|t_2|}$, the leaf labeled $x$ is a leaf of $t_2$ and not of $t_1$ and in this case $D_x(t)=t_1*D_x(t_2)$.

Note that, conditional on the set of leaf labels they have, $t_1$ and $t_2$ are each random trees which are invariant under any permutation of their respective leaf labels.
Thus, applying $F_l$ to $t$ and using the inductive assumption shows that $D(t) = \frac{|t_1|}{|t_1|+|t_2|}D(t_1)*t_2 + \frac{|t_2|}{|t_1|+|t_2|}t_1*D(t_2)$ as desired.

Combining these conditioned results over all possible first splits gives the desired inductive step.

The thin case now follows by commutativity of the forgetful functions \ignore{(Proposition \ref{prop:forget-commute})}, Proposition \ref{prop:delete-forget-orientation-commute} and the symmetry of Definition \ref{def:delete-x}.
\end{proof}

\begin{definition}
A sequence of probabilities $(P_n)_{n=0}^\infty$, such that $P_n$ is a probability on (fat or thin) cladograms with $n$ leaves, is called {\em deletion stable} if $D_n(P_n) = P_{n-1}$ for all $n\ge1$.
\end{definition}

\begin{definition}
A sequence of probabilities $(P_n)_{n=0}^\infty$, such that $P_n$ is a probability on (fat or thin) tree shapes with $n$ leaves, is called {\em deletion stable} if $D(P_n) = P_{n-1}$ for all $n\ge1$.
\end{definition}

\begin{corollary}
\label{cor:deletion-stable-forget-commutes}
If $(P_n)_{n=0}^\infty$ is a sequence, such that $P_n$ is a probability on (fat or thin) cladograms with $n$ leaves, which is {\em deletion stable} and invariant under permutation of leaf labels then the sequence $(F_l(P_n))_{n=0}^\infty$ of probabilities on (fat or thin) tree shapes is deletion stable.
\end{corollary}
\begin{proof}
This follows directly from the previous two definitions and Proposition \ref{prop:delete-forget-commute}.
\end{proof}

\subsection{Deletion stability and conditional split probabilities}
In the case of Markovian self-similar probabilities, deletion stability is equivalent to the conditional split distribution satisfying a simple `consistency condition'.
This condition is used in the next section to show that the alpha models are deletion stable.

Recall that if $q$ is a conditional split probability then it must satisfy $\sum_{m=1}^{n-1}q(m,n-m) = 1$ for all integers $n\ge2$.

\begin{proposition}
\label{prop:deletion-stability-if-mss}
Let $S=(P_n)_{n=0}^\infty$ be a sequence, such that $P_n$ a probability on (fat or thin) tree shapes with $n$ leaves, or a sequence such that $P_n$ is a probability on (fat or thin) cladograms which is invariant under permutations of leaf labels.
If $S$ is Markovian self-similar then it is deletion stable if and only if it has a conditional split distribution $q$ satisfying $q(x,y) =$
$$\frac{1}{1-\frac{q(1,x+y)+q(x+y,1)}{x+y+1}}\left(q(x+1,y) \frac{x+1}{x+y+1} + q(x,y+1)\frac{y+1}{x+y+1}\right)$$
for all integers $x,y \ge 1$.
\end{proposition}

\begin{proof}
First to reduce the cladogram cases to the tree shape cases.  If $S = (P_n)_{n=0}^{\infty}$ is a sequence on fat or thin cladograms then, by Corollary \ref{cor:deletion-stable-forget-commutes}, this sequence is deletion stable if and only if $(F_l(P_n))_{n=1}^\infty$ is deletion stable.  Since forgetting leaf labels leaves the conditional split distribution unchanged (Definition \ref{def:MSS-2}), proving the cladogram cases reduces to proving the cases of fat or thin tree shapes.

In the case of fat tree shapes, the conditional split distribution, $q$, of $S$ is uniquely defined (Proposition \ref{prop:unique-symmetric-q}).
In the case of thin tree shapes, there is a unique conditional split distribution $q$ of $S$ which is symmetric in the sense that $q(a,b) = q(b,a)$ (Proposition \ref{prop:unique-symmetric-q}).  Take this split distribution $q$.

By the definition of the conditional split distribution, the probability measure $P_{n+1} = \sum_{m=1}^{n}q(m,n+1-m)P_m*P_{n+1-m}$ for all $n\ge1$ and so
$$D(P_{n+1}) = \sum_{m=1}^{n}q(m,n+1-m) D(P_m*P_{n+1-m})$$
for all $n\ge1$.
By the definition of $D$ this is equal to:
$$\sum_{m=1}^{n}q(m,n+1-m) ( \frac{m}{n}D(P_m)*P_{n+1-m} + \frac{n+1-m}{n}P_m*D(P_{n+1-m}))$$

Noting that $P_k*P_0 = P_0*P_k$ for all $k$, this expression may be rearranged into
$$
P_{n}\frac{1}{n+1}(q(1,n)+q(n,1)) +$$ $$
\sum_{m=1}^{n} \left( q(m+1,n-m) \frac{m+1}{n+1}D(P_{m+1})*P_{n-m}
 + q(m,n-m+1)\frac{n-m+1}{n+1} P_m*D(P_{n-m+1}) \right)
 $$
Let $c_{n+1} = 1-\frac{1}{n+1}(q(1,n)+q(n,1))$.

Thus, if $S = (P_n)_{n=0}^\infty$ is deletion stable then since $D(P_{k+1}) = P_k$ for all $k\ge0$ it follows that:
$$P_n c_{n+1} = 
\sum_{m=1}^{n-1} \left(q(m+1,n-m) \frac{m+1}{n+1} + q(m,n-m+1)\frac{n-m+1}{n+1}\right) P_m*P_{n-m}$$

Note that if $q(a,b) = q(b,a)$ then the coefficients of $P_m*P_{n-m}$ and $P_{n-m}*P_m$ on the right hand side are equal.  

Thus, the uniqueness of the conditional split distribution $q$ in the case of fat tree shapes, and the uniqueness of the symmetric conditional split distribution in the case of thin tree shapes, implies that
$$q(m,n-m)c_{n+1} = q(m+1,n-m) \frac{m+1}{n+1} + q(m,n-m+1)\frac{n-m+1}{n+1}$$
for all integers $n,m \ge 1$ such that $m<n$.

On the other hand, suppose that $q$ is a conditional split distribution of $S = (P_n)_{n=0}^\infty$ which satisfies the equation in the statement of this proposition.  Induction shows that $(P_n)_{n=0}^\infty$ is deletion stable as follows: 

It is always true that $D(P_1) = P_0$ and $D(P_2) = P_1$ as there are unique (fat or thin) tree shapes with $0,1$ and $2$ leaves.  Suppose that $D(P_{k}) = D(P_{k-1})$ for all $k \in \{1,2,\ldots,n\}$.  Then, by the computations above,
$$D(P_{n+1}) = P_{n}\frac{1}{n+1}(q(1,n)+q(n,1)) +$$ $$
\sum_{m=1}^{n} \left( q(m+1,n-m) \frac{m+1}{n+1}D(P_{m+1})*P_{n-m}
 + q(m,n-m+1)\frac{n-m+1}{n+1} P_m*D(P_{n-m+1}) \right)
 $$
which by the inductive assumption is equal to:
\begin{eqnarray}
=  & \sum_{m=1}^{n} \left( \frac{1}{n+1}(q(1,n)+q(n,1))q(m,n-m) \right. \nonumber \\
& \left. + \frac{(m+1)q(m+1,n-m) + (n-m+1)q(m,n-m+1)}{n+1} \right) P_m*P_{n-m} \nonumber
\end{eqnarray}

By the definition of the conditional split distribution $q$ and the assumption that it satisfies the equations given in the statement of the proposition, this expression is equal to $P_n$.  Thus $D(P_{n+1}) = P_n$, and so by induction this holds for all $n>0$.
\end{proof}

\subsection{The case of the alpha model}
\begin{proposition}
All four of the alpha models are deletion stable for every value of alpha in $[0,1]$.
\end{proposition}
\begin{proof}
Recall that the alpha models on fat and thin cladograms are invariant under permutations of leaf labels.
Lemma \ref{lemma:alpha-mss} states that all four of the alpha models are Markovian self-similar and that the conditional split distributions, $q_\alpha$, for the alpha models are given by 
$$
q_\alpha(a,b) = \frac{\Gamma_\alpha(a)\Gamma_\alpha(b)}{\Gamma_\alpha(a+b)}
\left( \frac{\alpha}{2}\binom{a+b}{a} + (1-2\alpha)\binom{a+b-2}{a-1} \right)
$$
$$
= \frac{\Gamma_\alpha(a)\Gamma_\alpha(b)}{\Gamma_\alpha(a+b)}\binom{a+b}{a}
\left( \frac{\alpha}{2} + (1-2\alpha)\frac{ab}{(a+b)(a+b-1)} \right)
$$

It remains to show that $q_\alpha$ satisfies the equations given in Proposition \ref{prop:deletion-stability-if-mss}, for all values of $\alpha$ in $[0,1]$.

Let $a,b \ge 1$, and $n=a+b$.  Let  $A = q_\alpha(a+1,b) \frac{a+1}{a+b+1} + q_\alpha(a,b+1)\frac{b+1}{a+b+1}$.  It is sufficient to show that $A = \left(1- \frac{1}{n+1}(q_\alpha(1,n)+q_\alpha(n,1))\right) q_\alpha(a,b)$

Expanding A and rearranging gives:
$$A =
\frac{\Gamma_\alpha(a+1)\Gamma_\alpha(b)}{\Gamma_\alpha(a+b+1)}\left( \frac{\alpha}{2}\binom{a+b+1}{a+1} + (1-2\alpha)\binom{a+b+1-2}{a}\right)\frac{a+1}{a+b+1}$$
$$+ \frac{\Gamma_\alpha(a)\Gamma_\alpha(b+1)}{\Gamma_\alpha(a+b+1)}\left( \frac{\alpha}{2}\binom{a+b+1}{b+1} + (1-2\alpha)\binom{a+b+1-2}{b}\right)\frac{b+1}{a+b+1}$$

$$= \frac{\Gamma_\alpha(a)\Gamma_\alpha(b)}{\Gamma_\alpha(a+b+1)}\binom{a+b}{a}
\left( (a-\alpha)\frac{\alpha}{2} + (a-\alpha)(1-2\alpha)\frac{(a+1)b}{(a+b+1)(a+b)} \right)
$$
$$+\frac{\Gamma_\alpha(a)\Gamma_\alpha(b)}{\Gamma_\alpha(a+b+1)}\binom{a+b}{a}
\left( (b-\alpha)\frac{\alpha}{2} + (b-\alpha)(1-2\alpha)\frac{(b+1)a}{(a+b+1)(a+b)} \right)
$$

Let $C = \frac{\Gamma_\alpha(a)\Gamma_\alpha(b)}{\Gamma_\alpha(a+b+1)}\binom{a+b}{a}\frac{1}{a+b+1}$
So that
$$A = C \times \left( \frac{\alpha}{2}(a+b-2\alpha)(a+b+1) + \frac{(1-2\alpha)\left( (a-\alpha)(a+1)b + (b-\alpha)(b+1)a\right)}{a+b}\right)$$
which may be rearranged into:
$$ = C \times \frac{(a+b+2(1-\alpha))(\alpha(a+b)(a+b-1) +2ab)}{2(a+b)}$$

Now, $1- \frac{1}{n+1}(q(1,n)+q(n,1)) = 1 - \frac{1}{n+1}\frac{2}{n-\alpha}\left( \frac{\alpha}{2} + (1-2\alpha)\right)$.  Since $n=a+b$, this is equal to $\frac{(a+b-1)(a+b+2(1-\alpha))}{(a+b+1)(a+b-\alpha)}$.

Expanding $q(a,b)$ and rearranging gives:
$$q(a,b) = \frac{\Gamma_\alpha(a)\Gamma_\alpha(b)}{\Gamma_\alpha(a+b)} \left( \frac{\alpha}{2}\binom{a+b}{a} + (1-2\alpha)\binom{a+b-2}{a-1}\right)$$
$$= \frac{\Gamma_\alpha(a)\Gamma_\alpha(b)}{\Gamma_\alpha(a+b+1)} \binom{a+b}{a}\frac{1}{a+b+1}
\times (a+b-\alpha)(a+b+1)
\left( \frac{\alpha}{2} + (1-2\alpha)\frac{ab}{(a+b)(a+b-1)}\right) $$
$$= C \times \frac{(a+b-\alpha)(a+b+1)(\alpha(a+b)(a+b-1) +2ab)}{2(a+b)(a+b-1)}$$

Thus
$\left(1- \frac{1}{n+1}(q(1,n)+q(n,1))\right)q(a,b) = C \times \frac{(a+b+2(1-\alpha))(\alpha(a+b)(a+b-1) +2ab)}{2(a+b)} = A$
as desired.

Thus the conditions of Proposition \ref{prop:deletion-stability-if-mss} are satisfied for all four alpha models, and so they are deletion stable.
\end{proof}

Although perfectly correct, the above proof does not provide a good intuitive sense of {\em why} the alpha models are deletion stable.

One answer to this is to view the alpha model as the stationary distributions of the {\em delete-alpha-insert Markov chains}.  These will be discussed in a subsequent paper.
\ignore{
These are discussed in Section \ref{chapter:DAI-markov-chains}.  They fit together in such a way as to provide an easy proof that the stationary distribution of one is obtained from the others by inserting or deleting the correct number of leaves.
}
\subsection{A note on multifurcating trees}
The general definitions and results of this chapter may all be extended to multifurcating trees.
In particular, the definitions of Markovian self-similarity, conditional split distribution, deletion of a uniform random leaf or labeled leaf, and deletion stability all extend in the obvious way.
There is also a natural extension of Proposition \ref{prop:deletion-stability-if-mss} to the case of multifurcating trees.  The conditions on the split distribution are natural extensions of those for binary trees.

For the sake of brevity, this material is omitted.

\subsection{Other probabilities on Cladograms}
The alpha models are some of many different probabilities on cladograms and tree shapes.  The most popular and well known of these are the Yule, Uniform and Comb models.  These three are also Markovian self-similar and deletion stable.
The only other known models with these properties are the alpha model described here and the betal model of Aldous.

A major attraction of the alpha model is that it interpolates smoothly between the Yule, Uniform and Comb models.  The beta model of Aldous also interpolates between these three and extends beyond the Yule model to give models with very flat trees.

These models are now briefly discussed and compared with the alpha model.
\subsection{The Yule, Uniform and Comb models}
The {\em Yule model}, or {\em neutral evolution model}, was first defined by Yule in 1924 \cite{Yule-1924}.  It may be described in many different ways.  The most convenient description here is the following (see  \cite{Athreya-Ney-1972}):  Starting with a single species/leaf, at each step choose one of the extant species to bifurcate (split into two species) until the required number of species is reached.

The {\em Uniform model} is simply the uniform distribution on cladograms of a given size.  It is well known that there are $(2n-3)!!$ cladograms with exactly $n\ge2$ leaves.  See \cite{Cavalli-Sforza-Edwards-1967} for example.

The {\em Comb model} is the sequence of probabilities which assign probability $1$ to the most asymmetric tree of each size, called the {\em comb tree}.

From the above description of the Yule model, it is clear that this is precisely the alpha model with $\alpha=0$ since every new leaf is inserted at a uniform random leaf edge.  Similarly, a simple induction shows that setting $\alpha=1/2$ gives the Uniform model, as the next leaf is inserted at a uniformly chosen edge.
Finally, setting $\alpha=1$ gives the Comb model since every new leaf is inserted at a uniform random internal edge.

A more formal proof of this fact goes as follows:

When $\alpha=0$ the conditional split distribution of the alpha model is 
\begin{eqnarray}
q_0(a,b) &=& \frac{\Gamma_0(a)\Gamma_0(b)}{\Gamma_0(a+b)}
\left( \frac{0}{2}\binom{a+b}{a} + (1-2\times 0)\binom{a+b-2}{a-1} \right) \nonumber \\
&=& \frac{(a-1)! (b-1)!}{(a+b-1)!} \frac{(a+b-2)!}{(a-1)!(b-1)!} \nonumber \\
&=& \frac{1}{a+b-1} \nonumber
\end{eqnarray}
which is split distribution of the Yule model.

For the case of the Uniform model, a simple counting argument shows that the conditional split distribution satisfies $q(a,b) = \binom{a+b}{a} \frac{c_a c_b}{c_{a+b}}$, where $c_n = (2n-3)!!$ is the number of cladograms with $n$ leaves.

When $\alpha=\frac{1}{2}$ the conditional split distribution for the alpha model is:
\begin{eqnarray}
q_{\frac{1}{2}}(a,b) &=& \frac{\Gamma_{\frac{1}{2}}(a)\Gamma_{\frac{1}{2}}(b)}{\Gamma_{\frac{1}{2}}(a+b)}
\left( \frac{{\frac{1}{2}}}{2}\binom{a+b}{a} + (1-2\times {\frac{1}{2}})\binom{a+b-2}{a-1} \right) \nonumber \\
&=&\frac{(a-1-{\frac{1}{2}})\ldots (1-{\frac{1}{2}}) (b-1-{\frac{1}{2}})\ldots(1-{\frac{1}{2}})}{(a+b-1-{\frac{1}{2}})\ldots(1-{\frac{1}{2}})} \frac{1}{4} \binom{a+b}{a} \nonumber\\
&=&\frac{1}{2} \frac{(2a-3)(2a-5)\ldots(3)(1)(2b-3)\ldots(3)(1)}{(2(a+b)-3)(2(a+b)-5)\ldots(3)(1)} \nonumber\\
&=&\frac{1}{2} \binom{a+b}{a}\frac{(2a-3)!!(2b-3)!!}{(2(a+b)-3)!!} \nonumber\\
&=&\frac{1}{2} \binom{a+b}{a}\frac{c_a c_b}{c_{a+b}} \nonumber
\end{eqnarray}
which is the conditional split distribution of the Uniform model.

When $\alpha=1$ the conditional split distribution of the alpha model is $q_1(1,n) = q_1(n,1) = \frac{1}{2}$ for $n>1$ and $q_1(a,b)=0$ if $a$ or $b$ is not equal to $1$.  This is the conditional split distribution of the Comb model.

Notice that Yule trees tend to be flatter than Uniform trees, which are of course flatter than the Comb tree.  Similarly, the average depth of leaves in a Yule tree is less than that in a Uniform tree which is less than that in a Comb tree. These observations can be made more precise using Colless' and Sackin's index.  The inequalities extend to the alpha model and are made precise in Section \ref{chapter:Sackins-Colless-index}.

\subsection{The beta model of Aldous}
\label{sec:beta-model}
The other probabilities on cladograms of interest are those of the {\em beta model} of David Aldous, described in \cite{Aldous-1996}.  Other than the alpha models, the beta model is the only known family which interpolates between the Yule, Uniform and Comb models and is Markovian self-similar and deletion stable.
The beta model flows from a different description of the Yule model: uniform stick breaking.  This uniform stick breaking is extended to stick breaking according to the beta distribution on the unit interval.  The conditional split probabilities which arise are then extended beyond the point where stick breaking make sense.

Like the alpha model, the beta model is deletion stable (sampling consistent) and Markovian self-similar.  It is parameterized  by a single variable $\beta \in (-2,\infty]$, passes through the Yule model at $\beta=0$, the Uniform models at $\beta=-\frac{3}{2}$, and converges to the Comb model as $\beta \rightarrow -2$.  Unlike the alpha model, it extends beyond the Yule model to give distributions with much flatter trees ($\beta>0$).

As $\beta \rightarrow \infty$ it converges to the model defined by `perfect $1/2:1/2$ stick breaking'.  This model should be the `flattest possible' sampling consistent, Markovian self-similar distribution on cladograms.  Here `flattest possible' can mean either lowest expected value of Colless' (or Sackin's) index for all sizes of cladogram.

\subsection{The alpha model is not the beta model}
Here is a short proof that the alpha and beta models are different, and in fact only intersect at the Yule, Uniform and Comb models.

The conditional split distribution of the beta model is
$$q(a,b) = \frac{1}{k_n(\beta)}\frac{\Gamma(\beta +a+1)\Gamma(\beta+b+1)}{\Gamma(a+1)\Gamma(b+1)}$$
where $k_n(\beta)$ is a normalizing constant.  This is given in \cite{Aldous-1996} and \cite{Aldous-2001}.

\begin{theorem}
The alpha model and the beta model intersect only at the Yule, Uniform and Comb models.
\end{theorem}
\begin{proof}
Since it has already been shown that both models pass through the Yule, Uniform and Comb models, all that remains is to show that they do not intersect at any other point.  It is sufficient to show that at no other point do the conditional split distributions agree.

Consider the conditional split distribution for six leaves.  To avoid dealing with the normalization constant in the beta model, take the ratios $\frac{q(1,5)}{q(2,4)}$ and $\frac{q(2,4)}{q(3,3)}$.

For the alpha model these ratios are, respectively, $\frac{2(1+\alpha)(4-\alpha)}{(1-\alpha)(8-\alpha)}$ and $\frac{(8-\alpha)}{4(2-\alpha)}$.

For the beta model these ratios are, respectively, $\frac{\beta+5}{\beta+2}\frac{2}{5}$ and $\frac{\beta+4}{\beta+3}\frac{3}{4}$.

Equating the first ratio of split probabilities gives:
$$\frac{2(1+\alpha)(4-\alpha)}{(1-\alpha)(8-\alpha)}=\frac{\beta+5}{\beta+2}\frac{2}{5}$$
Solving for $\beta$ gives
$$\beta = \frac{5\alpha(5-\alpha)}{6(\alpha^2-4\alpha-2)}$$
Equating the second ratio of split probabilities gives:
$$\frac{(8-\alpha)}{4(2-\alpha)}=\frac{\beta+4}{\beta+3}\frac{3}{4}$$
Solving for $\beta$ gives:
$$\beta = \frac{-9\alpha}{2(1+\alpha)}$$
Thus, if the two models are equal it must be that:
$$\frac{5\alpha(5-\alpha)}{6(\alpha^2-4\alpha-2)}= \frac{-9\alpha}{2(1+\alpha)}$$
In other words:
$$\frac{ -2\alpha^3 + 8\alpha^2 - \frac{7}{2}\alpha}{\alpha^3 - 2\alpha^2 - 6\alpha -2} = 0$$
Which happens only if
$$2-\alpha^3 + 8\alpha^2 - \frac{7}{2}\alpha =0$$
Solving for $\alpha$ gives
$\alpha = 0,\frac{1}{2},\frac{7}{2}$.

Since the alpha model is not defined for $\alpha = \frac{7}{2}$ and the other two values correspond the the Yule and Uniform model this completes the proof.
As a final note, $\alpha=1$ did not appear as a solution because in that case (and only that case) the ratios are not real numbers.
\end{proof}

\section{Sackin's index and Colless' index}
\label{chapter:Sackins-Colless-index}
This section addresses two common statistics of tree shape: Sackin's index and Colless' index.
Sackin's index is the sum of the depth of all leaves in the tree.  In other words, the sum of the distance between the root and each leaf.
Colless' index is computed as follows:
For each internal vertex, compute the absolute value of the difference between the number of leaves below each of the two children, then sum up these numbers.

Sackin's index dates back to a paper of M.J. Sackin in 1972 \cite{Sackin-1972}, and Colless' to a paper of his in 1982 \cite{Colless-1982}.  These indices and others are described in an excellent paper of Shao and Sokal \cite{Shao-Sokal-1990}.
Formal symbolic definitions of each of these indices are given below.

In this chapter, the expected value of both of these indices for a cladogram of size $n$ chosen according to the alpha model is shown to be $O(n^{1+\alpha})$ for $\alpha \in (0,1]$ and $O(n \log n)$ for $\alpha=0$.  Dividing by $n$ shows that the expected depth of a random leaf is $O(n^\alpha)$.

Previous work on these and other indices in the cases of the Yule and Uniform models may be found in \cite{Shao-Sokal-1990} \cite{Heard-1992} \cite{Rogers-1993} \cite{Kirkpatrick-Slatkin-1993} \cite{Rogers-1996} \cite{Mooers-etal-1995} \cite{Mooers-Heard-1997} \cite{Blum-Francois-2005} \cite{Blum-Francois-Janson-2005}.

\subsection{Sackin's and Colless' indices defined}
Now for a formal definition of Sackin's index.
Denote by $S(T)$ the value of Sackin's index and $C(T)$ the value of Colless' index on a tree shape or cladogram $T$.  

Recall that the distance between two vertices in a tree is denoted by $d$.
\begin{definition}
Given a fat or thin tree shape $t$ with root $r$ and leaf set $s$, Sackin's index for this tree is defined to be $S(t) = \sum_{v \in s} (d(v,r) - 1)$.
\end{definition}
Note that the `depth' of a leaf is counted from the first branch point rather than the root: $d(r,v)-1$ rather than $d(r,v)$.  This is because many authors do not include the root edge in a tree shape, and also to preserve the alternative definition of Sackin's index given below.

For a vertex $v$ of a tree, let $N_v$ denote the number of leaves below and including $v$.
An equivalent definition of Sackin's index is:
\begin{definition}
Given a fat or thin tree shape $t$ with internal vertex set $I$, Sackin's index for this tree is defined to be $S(t) = \sum_{v\in I} N_v$.
\end{definition}

\begin{proposition}
The two preceding definitions of Sackin's index agree.
\end{proposition}
\begin{proof}
Let $t$ be a fat or thin tree shape with root $r$, leaf set $s$ and internal vertex set $I$.
Let $[S]$ denote the indicator function of a statement $S$.  In other words, $[S]=1$ if $S$ is true and $[S]=0$ otherwise.
Since $t$ is a tree, the path from a leaf to the root is unique and passes through every vertex above the leaf exactly once.  
Thus
$\sum_{v \in s} d(v,r) = \sum_{v \in s} \sum_{u \in I} [\text{$u$ above $v$}]$.  Exchanging the order of summation, this becomes $\sum_{u \in I} \sum_{v \in s} [\text{$u$ is above $v$}]$ which is equal to $\sum_{u \in I} \sum_{v \in s} [\text{$v$ is below $u$}]$ which by the definition of $N_u$ is equal to $\sum_{u \in I} N_u$.
Thus the two definitions of Sackin's index agree.
\end{proof}

Next to define Colless' index.  First some notation is introduced.
Every internal vertex, $v$, of a (fat or thin) tree shape has exactly two children.  Let $L_v$ denote the number of leaves below the left child and $R_v$ the number of leaves below the right child.  If the tree shape is thin then choose which child is 'left' and which is 'right' arbitrarily.
\begin{definition}
Given a (fat or thin) tree shape $t$ with internal vertex set $I$, Colless' index for this tree shape is defined to be $C(t) = \sum_{v\in I} |L_v - R_v|$
\end{definition}

Sackin's and Colless' indices for fat or thin cladograms are defined by first applying the map which forgets leaf labels and then calculating the index.

The following identity may be found in \cite{Blum-Francois-Janson-2005} and is used later in this chapter to show that Sackin's and Colless' indices have asymptotic covariance $1$ for all alpha models except $\alpha=0$.
\begin{lemma}
\label{lemma:Sackin-vs-Colless}
If $t$ is a fat or thin tree shape with internal vertex set $I$ then $C(t) = S(t) - 2 \sum_{v\in I} \min(L_v,R_v)$.
\end{lemma}
\begin{proof}
By the definition of Colless' index, $C(t) = \sum_{v\in I} |L_v - R_v| = \sum_{v\in I} L_v + R_v - 2 \min(L_v,R_v)$.  Since  $N_v = L_v + R_v$ for every internal vertex, $v$, it follows that $C(t) = \sum_{v\in I} N_v - 2 \min(L_v,R_v) = S(t) - 2 \sum_{v\in I} \min(L_v,R_v)$.
\end{proof}

Now, each of these two maps, $S(t)$ and $C(t)$, may be applied to a random variable on tree shapes to give a real random variable representing the distribution of each statistic.
Let $S_n(\alpha)$ and $C_n(\alpha)$ denote the random variables arising in this way from a random variable on tree shapes with $n$ leaves which is distributed according to the alpha model on trees with $n$ leaves.
In other words, for a tree with $n$ leaves chosen randomly under the alpha model, let $S_n(\alpha)$ denote the distribution of Sackin's index and let $C_n(\alpha)$ denote the distribution of Colless' index.

These random variables have already been studied in great detail in the cases of the Yule ($\alpha=0$) and Uniform ($\alpha=1/2$) models.  Results for these cases are surveyed in the next subsection.  Some of these results are then generalized to cover all values of alpha in $[0,1]$.

\subsection{Sackin's and Colless' index for alpha trees}
In this section, some of the results just quoted will be generalized to all values of alpha.
In particular, the expected value of Sackin's index for an alpha tree with $n$ leaves is $S_n(\alpha) = O(n^{1+\alpha})$ for $\alpha\in(0,1]$.
This implies that, for $\alpha\in (0,1]$, Colless' index is also $O(n^{1+\alpha})$ and the covariance of Sackin's index and Colless' index is $1$.

\subsection{The Yule and Uniform cases}
Much is already known about the distribution of Sackin's index and Colless' index in the cases of the Yule ($\alpha=0$) and Uniform ($\alpha=1/2$) models.  In particular, the mean, variance and covariance are known. In the case of the Uniform distribution the limiting distribution, after rescaling is the Airy distribution.
These results are summarized or proven in the preprints of Blum, Francois and Janson \cite{Blum-Francois-2005}, \cite{Blum-Francois-Janson-2005}.  Several papers have presented estimates of these values attained by simulation, such as those of Rogers \cite{Rogers-1993}, \cite{Rogers-1996}.

In the case of the Yule model ($\alpha=0$):
The correctly normalized Sackin's index, $\frac{S_n(0) - {\mathbb E}S_n(0)}{n}$, converges in distribution as $n$ approaches infinity.  The limiting distribution satisfies a fixed-point equation given by Rosler in \cite{Rosler-1992}, and has variance $\sigma = 7 - \frac{2\pi^2}{3}$.

In the case of Uniform trees ($\alpha = 1/2$):
$(\frac{S_n(1/2)}{n^{3/2}}, \frac{C_n(1/2)}{n^{3/2}})$ converges in distribution to $(A,A)$, where $A$ is the Airy distribution.
This is proven in \cite{Blum-Francois-Janson-2005}.
It also follows directly from the work of Aldous on {\em continuum random trees}: \cite{Aldous-CRT1}, \cite{Aldous-CRT2}, \cite{Aldous-CRT3}.

Notice that  the mean and variance of $S_n(1/2)$ and $C_n(1/2)$ are both order $n^{1+1/2}$ and their covariance trends to $1$.

\subsection{The expected value of Sackin's index}
Now to show that the expected value of Sackin's index is $O(n^{1+\alpha})$ for $\alpha \in (0,1]$.
Begin by defining some new statistics on trees which are close to Sackin's index.  Next, find a recurrence equation which they satisfy and try to solve it.

For a tree shape or cladogram $t$ define:

\noindent $T(t) = \text{sum of leaf depths}$

\noindent $K(t) = \text{sum of depths of all internal vertices}$

\noindent $L(t) = \text{sum of the number of internal nodes below and including each internal node}$

Here the `depth' of a vertex is the number of edges in the shortest path between it and the root vertex: $d(r,v)$.

Let $T_\alpha(n)$, $L_\alpha(n)$, and $K_\alpha(n)$ denote the expectations of these variables under the alpha model on tree shapes with $n$ leaves.

Note that, in the notation of the previous section, $L(t) =  \sum_{v\in I} (N_v-1)$ since the number of internal nodes below a vertex is one less than the number of leaves for a binary tree.  Thus $L(t)$ for a tree $t$ is Sackin's index minus $n-1$, the number of internal vertices.

The first few values for each of these functions are:

\begin{tabular}{llll}
n & $T_\alpha(n)$ & $K_\alpha(n)$ & $L_\alpha(n)$ \\
1 & 1 & 0 & 0 \\
2 & 4 & 1 & 1 \\
3 & 8 & 3 & 3 \\
4 & $12 + \frac{2}{3-\alpha}$ & $5 + \frac{2}{3-\alpha}$ & $ 5 + \frac{2}{3-\alpha}$
\end{tabular}

Notice that for these small values $K_\alpha(n) = L_\alpha(n)$ and $T_\alpha(n) - K_\alpha(n) = 2n-1$.  
In fact, these relations hold for each individual tree:

\begin{proposition}
\label{prop:TKL}
For any binary rooted tree, $t$, with $n$ leaves:
\begin{itemize}
\item $T(t)$ is Sackin's index plus $n$
\item $T(t) - K(t) = 2n-1$, and 
\item $K(t) = L(t)$.
\end{itemize}
\end{proposition}

\begin{proof}
Let $s$ be the leaf set of tree $t$, $I$ the set of internal vertices and $r$ the root.
Now $T(t)$ is the sum of the distance from the root to each leaf, $\sum_{v\in s}d(r,s)$, and there are $n$ leaves. Thus it is equal to $\sum_{v\in s}(d(r,s)-1) + n$ which is Sackin's index plus $n$.

For the difference between the sum of leaf depths and the sum of internal node depths:
This is true for $n=1$.  Suppose that it is true for all trees with less than $k$ leaves.
Given a tree with $k$ leaves, the first split has $p$ leaves to the 'left' and $q$ leaves to the 'right' ($p+q=k$).
The difference between the total leaf depth and internal node depths on the left is $(2p-1) + p - (p-1) = 2p$, the difference on the right is $(2q-1) + q - (q-1) = 2q$.  Adding these together and subtracting $1$ for the depth of the first branching node gives a total difference of $2p+2q-1 = 2k-1$.

For the second part, if $J$ is the set of internal nodes then let $b(i,j)$ be $1$ if $i$ is above $j$ (closer to the root) and $0$ otherwise.  Then
$$\sum_{i \in J} \left({\sum_{j \in J} b(i,j) }\right) = \sum_{j \in J} \left({\sum_{i \in J} b(i,j) }\right)$$
The left hand side of this equation is the sum of the number of internal nodes below and including the internal node, and the right hand side is the sum of the depths of each internal node.
\end{proof}

A recurrence relation for the expected value of $L(t)$ under the alpha model is now derived.

\begin{proposition}
\begin{equation}
\label{eqn:L-recurrence}
L_\alpha(n+1)= L_\alpha(n)\frac{n+1}{n-\alpha} + \frac{(2n-1)(1-\alpha)}{n-\alpha}
\end{equation}
\end{proposition}
\begin{proof}
$L_\alpha$ satisfies the recurrence relation:

$L_\alpha(n+1) 
= \frac{n(1-\alpha)}{n-\alpha} \left({ L_\alpha(n) + 1 + \frac{L_\alpha(n) + (n-1)}{n} }\right)$

$ + \frac{(n-1)\alpha}{n-\alpha} \left({ L_\alpha(n) + (\frac{L_\alpha(n)}{n-1} -1)+ (\frac{L_\alpha(n)}{n-1} + 1) }\right)$

$= L_\alpha(n) + \frac{(n-1)\alpha}{n-\alpha} + \frac{1-\alpha}{n-\alpha}\left({ L_\alpha(n) + (n-1) }\right) + \frac{\alpha}{n-\alpha}\left({ 2L_\alpha(n) }\right)$

$= L_\alpha(n)\frac{n+1}{n-\alpha} + \frac{(2n-1)(1-\alpha)}{n-\alpha}$

The first collection of terms in the group corresponds to a leaf displacement, which occurs with probability $\frac{n(1-\alpha)}{n-\alpha}$.  When this occurs, all the old nodes are still above the nodes they were above before, contributing $L_\alpha(n)$.  The new internal node has exactly itself below itself and thus contributes $1$.  An existing internal node gains this new internal node as a descendant if it is above the displaced leaf, so this contribution is the equal to the expected number of internal nodes which are above the displaced leaf.  This is equal to the expected depth of the leaf minus $1$, which is $\frac{T_\alpha (n)}{n} - 1 = \frac{L_\alpha(n) + (n-1)}{n}$.

The second collection of terms corresponds to an internal node being displaced, and occurs with probability $\frac{(n-1)\alpha}{n-\alpha}$.  In this case, all of the old nodes are still above the nodes they were above before, contributing $L_\alpha(n)$.  The number of internal nodes the new node is below is equal to the number of nodes the one it displaced was below (excepting that node), for a total expected contribution of $\frac{L_\alpha(n)}{n-1} - 1$. The new internal node is above all the internal nodes that the number of the node it displaced was above plus one for itself (for an expected contribution of $\frac{L_\alpha(n-1)}{n-1}+1$).
\end{proof}

Similar arguments give recurrences for $T_\alpha$ and $K_\alpha$.  The resulting recurrences may also be obtained by substituting the equations in Proposition \ref{prop:TKL} into the equation in Equation \ref{eqn:L-recurrence}

Notice that $L_\alpha(n)$ is strictly increasing in $n$ for fixed $\alpha$.  This follows as $\frac{n+1}{n-\alpha}>1$ and $\frac{(2n-1)(1-\alpha)}{n-\alpha}>0$.  Also, $L_\alpha(n)$ is an increasing function of $\alpha$ for each fixed $n$, strictly increasing for $n\ge 4$.
To see this, note that both $\frac{1}{n-\alpha}$ both $\frac{(2n-1)(1-\alpha)}{n-\alpha}>0$ are strictly increasing in $\alpha$;  thus if $L_\alpha(n-1)$ is an increasing function of $\alpha$ then so is $L_\alpha(n)$.

\begin{theorem}
\label{thrm:Lrecurrence}
$L_\alpha(n)$ is $O(n^{1+\alpha})$ for $\alpha \in (0,1] $
\end{theorem}
\begin{proof}
Fix $\alpha \in (0,1]$.
Begin by showing that $L_\alpha(n)$ is $o(n^{1+\alpha+\epsilon})$ for all $\epsilon>0$.

Let $M(n) = L_\alpha(n) (1-\alpha)$ and
suppose that for some choice of $c$ and some sufficiently large $n$ it is true that $M(n) \le c n^{1+\alpha+\epsilon}$.
Then equation \ref{eqn:L-recurrence} gives:

$$M(n) \le c (n-1)^{1+\alpha+\epsilon} \frac{n+1}{n-\alpha} + \frac{2n-1}{n-\alpha}$$

letting $x=\frac{1}{n}$, this is

$$ \le c n^{1+\alpha+\epsilon} (1-x)^{1+\alpha+\epsilon}\frac{1+x}{1-\alpha x} + \frac{2-x}{1-\alpha x}$$

$$ \le c n^{1+\alpha+\epsilon} (1 - \epsilon x + O(x^2) )  + (2+(2-\alpha)x + o(x) )$$

$$ = c ( n^{1+\alpha+\epsilon} - \epsilon n^{\alpha+\epsilon} + o(n^{\alpha+\epsilon})) + 2+(2-\alpha)\frac{1}{n} + o(\frac{1}{n}) $$

For sufficiently large $n$ this gives:
$$ M(n) \le c n^{1+\alpha+\epsilon}$$

Applying induction starting at this value of $n$ gives the desired result.

Let $A(n) = M(n)/n^{1+\alpha}$. Thus $A(n)$ is $o(n^\epsilon)$ for all $\epsilon >0$.
As before, equation \ref{eqn:L-recurrence} leads to:
$$n^{1+\alpha} A(n) = n^{1+\alpha}(1+O(\frac{1}{n^2}))A(n-1) + 2 + o(1)$$
which rearranges into:
$$A(n) - A(n-1) = O(\frac{1}{n^2})A(n-1) + O(\frac{1}{n^{1+\alpha}})$$
This implies that $A(n)$ is bounded, as $\sum_{i=0}^\infty \frac{1}{i^k}$ is bounded for all $k>1$, in particular for $k = 1+\alpha$ and $k = 2-\epsilon$ (for sufficiently small $\epsilon$).  On the other hand, $A(n)$ is positive and is strictly increasing for sufficiently large $n$ and so is bounded away from $0$ by a definite amount.  Thus $A(n)$ is $O(1)$ and so $M(n)$ and $L(n)$ are $O(n^{1+\alpha})$.
\end{proof}

This immediately gives:
\begin{corollary}
\label{cor:alpha-Sakins-index}
For $\alpha \in (0,1]$ the expected value of Sakin's index for a random alpha tree with $n$ leaves is order $n^{1+\alpha}$.
\end{corollary}

Dividing by the total number of leaves, $n$, gives:
\begin{corollary}
\label{corr:alpha-average-leaf-depth}
For $\alpha \in (0,1]$, the expected depth of a random leaf in a tree chosen from the alpha model with $n$ leaves is $O(n^\alpha)$.
\end{corollary}

\subsection{Covariance of Sackin's and Colless' Index}
It has just been shown that the mean of $S_n(\alpha)$ is of order $n^{1+\alpha}$ for all $\alpha \in (0,1]$.  It will now be shown that in fact $\frac{S_n(\alpha) - C_n(\alpha)}{n^{1+\alpha}}$ converges to $0$ in probability for all $\alpha \in (0,1]$.

In fact, the values of Sackin's index and Colless' index on any tree of size $n$ differ by at most $n \log_2 n$.
This shortcuts the need for Lemma 3 in \cite{Blum-Francois-Janson-2005}, replacing it with an easier and much better result.

Given a tree shape or cladogram, $T$, define $v(T)$ to be the sum over all internal vertices of the minimum of the number of leaves below each child of that vertex. In other words, $v(T) = \sum_{v\in I} \min(L_v,R_v)$.  By Lemma \ref{lemma:Sackin-vs-Colless} this is half the difference between Sackin's index and Colless' index for the tree $T$.

It seems plausible for $v(T)$ to take its maximum value over all trees with a fixed number of leaves at a very `balanced' tree.  A perfectly balanced tree with $n=2^k$ leaves has value $v(T) = 2^{k-1}k = \frac{n \log_2 n}{2}$.  It also seems reasonable for this tree to have the greatest value over all trees with at most $2^k$ leaves.  This suggest that if $T$ is a tree with $n$ leaves then the difference between Colless and Sackin's index, $2v(T)$, is at most $n \log_2 n$.

This turns out to be a good heuristic.
\begin{lemma}
If $T$ is a tree shape or cladogram with $n$ leaves then the difference between Colless' and Sackin's index for $T$ is at most $n \log_2 n$.  Specifically $0 \le S(T) - C(T) \le n \log_2 n$. 
\end{lemma}
\begin{proof}
Recall that the difference between Sakin's and Colless' index for a tree $T$ is $S(T) - C(T) = 2\sum_{v\in I} \min(L_v,R_v) = 2v(T)$ (by Lemma \ref{lemma:Sackin-vs-Colless}), where $I$ is the set of internal nodes of $T$.
Let $f(n)$ be the maximum value of $v(T)$, over all tree shapes (or cladograms), $T$, with $n$ leaves.
Clearly $v(T) \ge 0$, and so $f(n) \ge 0$.  Now show that $f(n) \le \frac{n}{2} \log_2 n$.

The proof is by induction.  First, not that for $k=1,2,3$ there is only one tree shape with $k$ leaves and $v(T) = 0,1,2$ respectively.  These values are less than or equal to $(1 \log_2 1 )/2, (2\log_2 2)/2, (3\log_2 3)/2$ respectively.  Thus $f(k) \le (k \log_2 k)/2$ for $k=1,2,3$.

Suppose that $f(k) \le (k \log_2 k)/2$ for all $k \le n$.  Note that $f(n)$ satisfies the recurrence relation $f(n) = \max_{i \in \{1,2,3,\ldots,\lfloor n/2 \rfloor+1\} } f(i) + f(n-i) + i$.  This follows as for every tree $T$ with first split $\{i,n-i\}$ has $v(T)$ equal to $\min(i,n-i)$ plus the value of the left and right subtrees, which are bounded above by $f(i)$ and $f(n-i)$ respectively.  On the other hand, $f(i) + f(n-i) + \min(i,n-i)$ is obtained for the tree which is the root join of trees with $i$ leaves and $n-i$ which maximize $v$ for these numbers of leaves.

Assume without loss of generality that $i \le \frac{n}{2}$.
Thus it is sufficient to show that for all $1\le i \le \frac{n}{2}$, the following inequality holds:
$(n \log_2 n)/2 \ge ((n-i) \log_2 (n-i))/2 + (i \log_2 i)/2 + i$.
In other words, show that $0 \ge (n-i) \log_2 (n-i) + i \log_2 i  - n \log_2 n + 2i$.

The second derivative of the right hand side is $\frac{1}{i} + \frac{1}{n-i}$, which is always greater than zero. Thus the function is convex.  The inequality is true when $i=1$, and equality holds when $i=n/2$.  Therefor, by convexity, the equality holds for all $i$ between $1$ and $n/2$.

Thus the inductive step holds and the lemma is proven.
\end{proof}
In other words, for large trees which are not too symmetrical these two statistics are almost identical.

This leads immediately to the desired result:
\begin{corollary}
$\frac{S_n(\alpha) - C_n(\alpha)}{n^{1+\alpha}}$ converges to $0$ uniformly (and so in probability) as $n$ approaches $\infty$, for all $\alpha \in (0,1]$.
\end{corollary}

And also:
\begin{corollary}
$C_n(\alpha) = O(n^{1+\alpha})$, for all $\alpha \in (0,1]$.
\end{corollary}
\begin{proof}
This follows directly from the previous Corollary and Corollary \ref{cor:alpha-Sakins-index}.
\end{proof}

\section{Sweet Cherries}
\label{chapter:cherries}
An easily computed statistic of a cladogram is the number of {\em cherries}.  A {\em cherry} is a pair of leaves which are both adjacent to the same internal vertex.  For example, the balanced rooted tree with $4$ leaves has two cherries.

\epsfysize=2 cm			
\begin{figure}[ht]
\center{
	\leavevmode
	\epsfbox{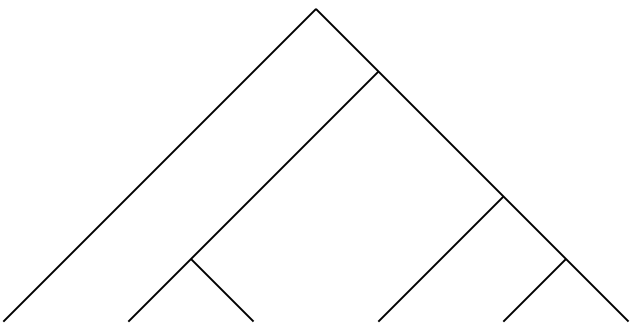}
	\caption{A tree with two cherries}
}
\end{figure}

McKenzie and Steel \cite{McKenzie-Steel-2000} showed that for the Yule model on rooted trees and Uniform model on unrooted trees the number of cherries is asymptotically normal, with known mean and variance.  These results are now extended to the alpha model:
\begin{theorem}
\label{thrm:cherries-normal}
If $\cherries_m$ is the number of cherries in a random alpha tree with $m$ leaves then for $\alpha \in [0,1)$
$$\frac{\cherries_m - m\frac{1-\alpha}{3-2\alpha}}{\sqrt{m\frac{(1-\alpha)(2-\alpha)}{(3-2\alpha)^2(5-4\alpha)}}} \longrightarrow {\cal N}(0,1)$$
\end{theorem}
For $\alpha=1$ and $m\ge2$, $\cherries_m$ is identically $1$ as a comb tree has only one cherry.

The proof of Theorem \ref{thrm:cherries-normal} follows the methods in \cite{McKenzie-Steel-2000}.
First, describe the formation of cherries in terms of an {\em extended Polya Urn model} and apply a theorem which proves asymptotic normality.
Next, use probability generating functions to find recurrences for the mean and variance.  Finally, solve these recurrences to find the asymptotic mean and variance.

Along the way, an exact formula for the mean is obtained.
The exact mean and variance have previously been calculated for the Yule model ($\alpha=0$) in \cite{McKenzie-Steel-2000}, and for the Uniform model on unrooted binary trees in \cite{Steel-Penny}, \cite{Hendy-Penny-1982} and \cite{Steel-1988} with these results collected in \cite{McKenzie-Steel-2000}.
An exact formula for the variance for all $\alpha$ may also be possible using the usual techniques for solving recurrence equations.

\subsection{Extended Polya urn models}
This section reviews a recent central limit theorem on {\em extended Polya urn} (EPU) models.  This result is to prove the asymptotic normality of the number of cherries in a random alpha tree.

First define the urn models.

Suppose an urn contains $k$ different types of balls.  If a ball of the $i$-th type is drawn from the urn then it is returned, along with $A_{ij}$ balls of the $j$-th type.    The value $A_{ij}$ may be negative, corresponding to the removal of balls from the urn.  Models with $A_{ii}>0$ are referred to as {\em generalized Polya urn} (GPU) models \cite{Athreya-Karlin}, \cite{Athreya-Ney-1972}.  Allowing for $A_{ii}$ to be negative, but requiring that the number of balls returned each time be a positive constant defines the class of {\em extended Polya urn} (EPU) models \cite{Bagchi-Pal}, \cite{Smythe}.

For both of these classes of urn models a number of asymptotic normality results exist.  The one relevant here (found in \cite{Bagchi-Pal}, \cite{Smythe}) is as follows:

\begin{theorem}
\label{thrm:EPU}
\cite{Bagchi-Pal} \cite{Smythe} Let $A = [A_{ij}]$ be the generating matrix for an EPU model, with principal eigenvalue $\lambda_1$.  Let $v$ be the left eigenvector of $A$ corresponding to $\lambda_1$, where the entries $v_i$ add up to one.  Also let $Z_{in}$ denote the number of balls of type $i$ in the urn after $n$ draws, where $i=1,2,\ldots,k$. For $k=2$ suppose that:

$(i)$ $A$ has constant row sums, where the constant is positive,

$(ii)$ $\lambda_1$ is positive, simple, and has a strictly positive left eigenvector $v$,

$(iii)$ $2\lambda < \lambda_1$ for the non-principal eigenvalue $\lambda$;

then $n^{-1/2}(Z_{1n} - n \lambda_1 v_1)$ has asymptotically a normal distribution with mean zero.

Furthermore, for $k>2$, suppose in addition:

$(iv)$ $2 \text{Re}(\lambda) < \lambda_1$ for all non-principal eigenvalues $\lambda$,

$(v)$ all complex eigenvalues are simple, and no two distinct complex eigenvalues have the same real part, except for conjugate pairs,

$(vi)$ all eigenvectors are linearly independent;

then $n^{-1/2}(Z_{1n} - n \lambda_1 v_1, Z_{2n} - n \lambda_1 v_2, \ldots, Z_{(k-1)n} - n \lambda_1 v_{(k-1)})$ has asymptotically a joint normal distribution with mean zero.
\end{theorem}

This theorem also applies to the case when the number of balls is a non-negative real number rather than a non-negative integer.

\subsection{The number of cherries is asymptotically normal}
This section follows the approach in \cite{McKenzie-Steel-2000} of describing the process of cherry formation in terms of an extended Polya urn model and applying theorem Theorem \ref{thrm:EPU}.

The asymptotic distribution for the number of cherries my be found by realizing the process of cherry formation as an EPU model.  Each new leaf is added in the alpha model by choosing an edge at random according to weights, breaking the edge in two with a new internal vertex and attaching a new leaf edge at that new vertex.  An extra cherry is created exactly when a leaf edge which is not already part of a cherry is chosen at the point of insertion.

\begin{proposition}
\label{prop:cherry-normal}
If $\cherries_m$ is the number of cherries in a random alpha tree with $m$ leaves then for $\alpha \in [0,1)$ there exists a variance $\sigma_n^2$ such that
$$\frac{\cherries_m - m\frac{1-\alpha}{3-2\alpha}}{\sigma_n} \rightarrow {\cal N}(0,1)$$
\end{proposition}
\begin{proof}
First to realize the creation of cherries as an extended Polya urn.

Let the first type of ball represent leaf edges which are part of a cherry, the second type of ball represent leaf edges which are not part of a cherry and the third type of ball represent internal edges.
Each non-cherry leaf edge is represented by a ball of type $2$ with weight $1-\alpha$ and each internal edge by a ball of type $3$ with weight $\alpha$.  Each cherry is represented by two balls of type $1$ with total weight $2(1-\alpha)$, as it consists of two leaf edges.

In this way, the total weight of all balls of a given type is proportional to the probability that the next leaf is inserted into that type of edge.
Note that the number of cherries is the weight of the first type of ball divided by $2(1-\alpha)$.

Now to determine what happens when a ball is chosen.

When a new leaf edge is inserted at a leaf edge which is already part of a cherry, the net effect is to add a new internal edge and a new non-cherry leaf edge.
The same happens when a new leaf edge is inserted at an internal edge.
When a new leaf edge is inserted at a leaf edge which is not part of a cherry then a new cherry is created, a new internal edge created, and a non-cherry leaf edge lost.  See Figure \ref{fig:cherry-formation}.

\begin{figure}[ht]
\label{fig:cherry-formation}
    \begin{center}
    \resizebox{7cm}{!}{\includegraphics{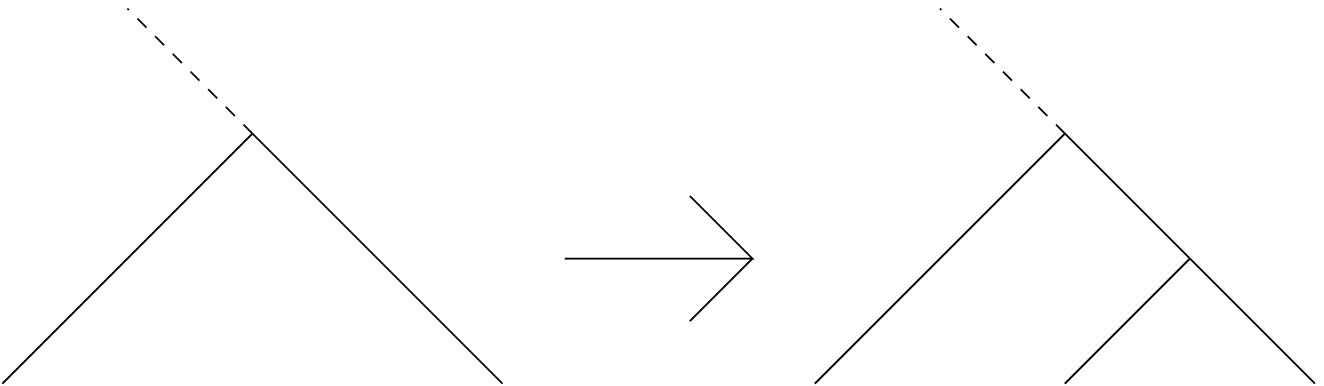}}\\
    \resizebox{8cm}{!}{\includegraphics{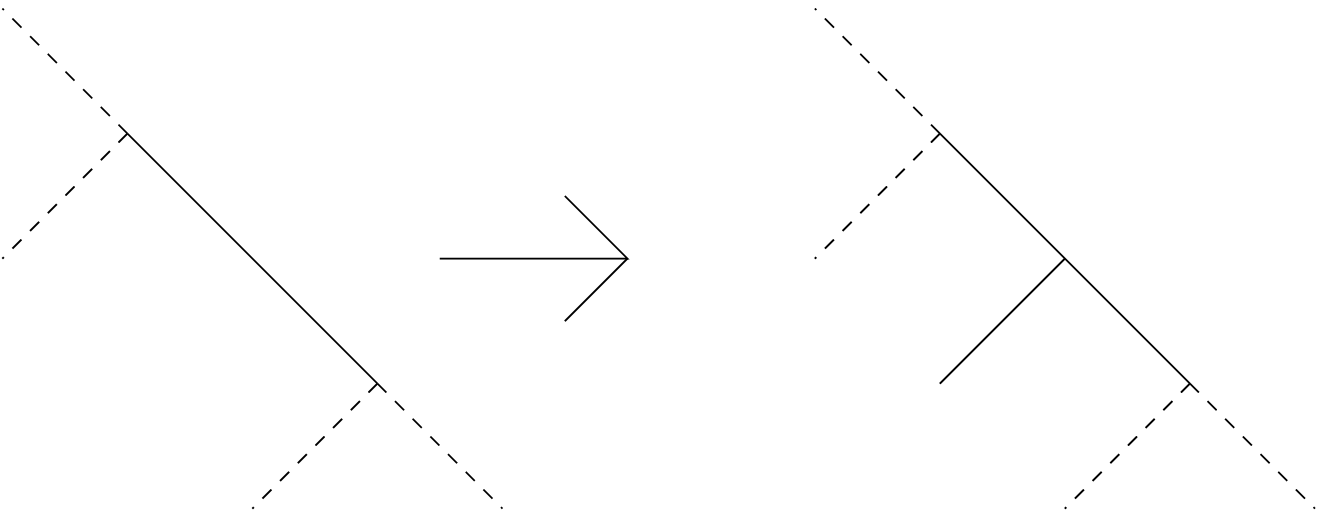}}\\
    \resizebox{7cm}{!}{\includegraphics{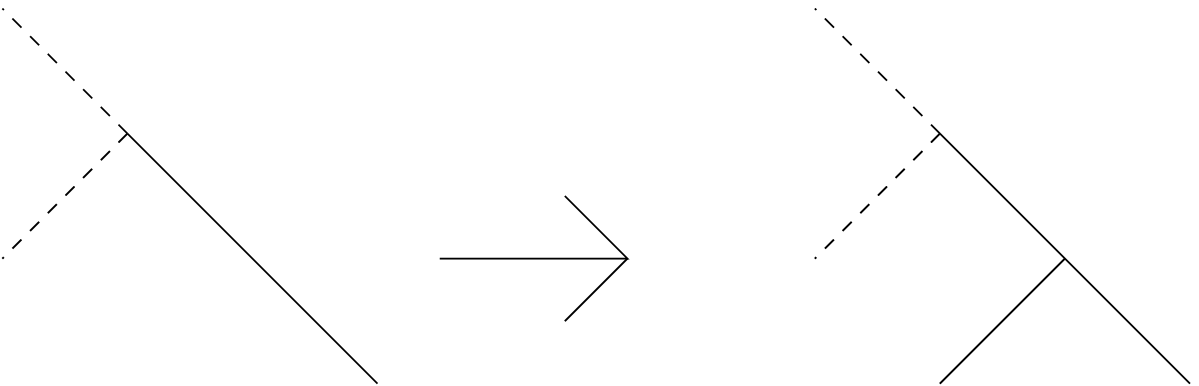}}\\	
	\caption[Cherry formation]{The effect of adding a leaf at a cherry leaf edge, an internal edge and a non-cherry leaf edge}
 \end{center}
\end{figure}

Recalling the weights chosen above, this means that the generating matrix for this urn scheme is:

$
A = \left[{
\begin{matrix}
0 & 1-\alpha & \alpha \\
2-2\alpha& -(1-\alpha) & \alpha \\
0 & 1-\alpha & \alpha
\end{matrix}
}\right]
$

This matrix has eigenvalues $1$, $0$ and $-2(1-\alpha)$, with corresponding eigenvectors $[\frac{2(1-\alpha)^2}{3-2\alpha},\frac{1-\alpha}{3-2\alpha},\alpha]$, $[1,0,-1]$ and $[1,-1,0]$.  As $\alpha \in [0,1]$ the principal eigenvalue is $\lambda_1 = 1$ and the corresponding left eigenvector, scaled such that its entries sum to one, is $[\frac{2(1-\alpha)^2}{3-2\alpha},\frac{1-\alpha}{3-2\alpha},\alpha]$.

Thus the conditions of the EPU asymptotics theorem, Theorem \ref{thrm:EPU}, are satisfied.
So, for some function $c$ of $m$,
$$\frac{1}{\sqrt{m}} \left({ Z_{1m} - m\frac{2(1-\alpha)^2}{3-2\alpha} }\right) \rightarrow {\cal N}(0,c)$$
where ${\cal N}(\mu,\sigma^2)$ is a normal distribution with mean $\mu$ and variance $\sigma^2$.

Finally, recall that $Z_{1m}$ is $2(1-\alpha)$ times the total number of cherries.
Therefore the desired result follows.
\end{proof}

\subsection{The mean and variance}
Recurrence equations are now found for the mean and variance of the number of cherries under the alpha model.
An exact formula for the mean is then found.

Let $\cherries_m$ be the number of cherries in a random tree shape or cladogram with $m$ leaves picked according to the alpha model.  Let $\mu_m$ be the mean of $\cherries_m$ and $\sigma_m^2$ the variance.  Note that each of these depends on the value of $\alpha$.
\begin{theorem}
\label{thrm:cherry_mean_variance}
The following recurrences hold:

$$\mu_{m+1} = \frac{m(1-\alpha)}{m-\alpha} + \frac{m-2+\alpha}{m-\alpha}\mu_m$$
$$
\begin{array}{rcl}
\sigma_{m+1}^2 & = & \frac{\alpha(1-\alpha)m(m-1)}{(m-\alpha)^2} + \sigma_m^2 \left({ \frac{m-4+3\alpha}{m-\alpha} }\right) \\
&& + \mu_m \left({ \frac{2(1-\alpha)(m(1-2\alpha) -\alpha)}{(m-\alpha)^2} }\right) - \mu_m^2 \frac{4(1-\alpha)^2}{(m-\alpha)^2}
\end{array}
$$
Furthermore:
\begin{equation}
\mu_m \sim m\frac{(1-\alpha)}{3-2\alpha}
;\hspace{10pt}
\sigma_m^2 \sim m\frac{(1-\alpha)(2-\alpha)}{(3-2\alpha)^2(5-4\alpha)}
\label{eqn:cherry_mean_variance}
\end{equation}
\end{theorem}
This theorem agrees with the corresponding theorems in \cite{McKenzie-Steel-2000} for the Yule ($\alpha=0$) and unrooted Uniform models (setting $\alpha=1/2$).

\begin{proof}
 
When inserting a new leaf into a cladogram the number of cherries increases if and only if the new leaf displaces a leaf which is not already part of a cherry.  If the number of cherries increases then it increases by exactly one.  Thus the variables $\cherries_m$ obey the following recurrence:

$$
\begin{array}{rcl}
\mathbb{P}[\cherries_{m+1}=k] &=& \mathbb{P}[\cherries_m=k-1] \frac{(1-\alpha)(m-2(k-1))}{m-\alpha} \\
& & + \mathbb{P}[\cherries_m=k] \frac{(m-1)\alpha + 2k(1-\alpha)}{m-\alpha}
\end{array}
$$

Let $P_m(x) = \sum_{k \ge 0} \mathbb{P}[\cherries_m=k]x^k$ be the probability generating function for $\cherries_m$.  Thus $P_1(x) = 1$ as the tree with $1$ leaf has no cherries, and $P_2(x)=x$ as the two leaf tree has exactly one cherry.

Now find a recurrence equation for $P_n(x)$.  The contribution to $P_{m+1}(x)$ from the first term in the above recurrence 
is:
$$\frac{m(1-\alpha)}{m-\alpha} x P_m(x) - \frac{2(1-\alpha)}{m-\alpha} x^2 \frac{d}{dx}P_m(x)$$
The contribution from the second term is:
$$\frac{(m-1)\alpha}{m-\alpha} P_m(x) + \frac{2(1-\alpha)}{m-\alpha} x \frac{d}{dx}P_m(x)$$

Thus the probability generating functions $P_m(x)$ satisfy the following recurrence equation:
\begin{equation}
P_{m+1}(x) = \frac{(m-1)\alpha + m(1-\alpha)x}{m-\alpha} P_m(x) + \frac{2(1-\alpha)}{m-\alpha}x(1 - x) \frac{d}{dx}P_m(x)
\label{eqn:cherries-probgenfn-recurrence}
\end{equation}

Note that $\mu_m = \frac{d}{dx}P_m(x)|_{x=1}$ and $\sigma_m^2 = \frac{d^2}{dx^2}P_m(x)|_{x=1} + \mu_m - \mu_m^2$.
For notational convenience let $P^{(k)}_m(x)$ denote $\frac{d^k}{dx^k}P_m(x)$.

Differentiating equation (\ref{eqn:cherries-probgenfn-recurrence}) yields:
$$
\begin{array}{rcl}
P_{m+1}^{(1)}(x) & = & \frac{m(1-\alpha)}{m-\alpha}P_m(x) + \frac{(m-1)\alpha + m(1-\alpha)x}{m-\alpha} P_m^{(1)}(x) \\
&& + \frac{2(1-\alpha)}{m-\alpha}(1-2x)P_m^{(1)}(x) + \frac{2(1-\alpha)}{m-\alpha}x(1-x)P_m^{(2)}(x)
\end{array}
$$

Evaluating at $x=1$, and noting that $P_m(1) = 1$ for all $m$, gives:

\begin{equation}
\label{eqn:cherry_mean_recurrence}
\mu_{m+1} = \frac{m(1-\alpha)}{m-\alpha} + \frac{m-2+\alpha}{m-\alpha}\mu_m
\end{equation}

There is one tree with two leaves and it has one cherry so $\mu_2 =1$.  By Proposition \ref{prop:cherry-normal}, that $\mu_m \sim \frac{1-\alpha}{3-2\alpha}$ so a direct solution is not presented here.

Differentiating equation (\ref{eqn:cherries-probgenfn-recurrence}) a second time gives:
\begin{equation}
\label{eqn:cherry_probgenfn_second_derivative}
\begin{array}{rcl}
P_{m+1}^{(2)}(x) & = & P_m^{(1)}(x) \frac{m(1-\alpha)}{m-\alpha} + P_m^{(1)}(x) \frac{m(1-\alpha) - 4(1-\alpha)}{m-\alpha}  \\
&& + P_m^{(2)}(x) \frac{(m-1)\alpha + m(1-\alpha)x + 2(1-\alpha)(1-2x)}{m-\alpha}
+ P_m^{(2)}(x) \frac{2(1-\alpha)}{m-\alpha}(1-2x) \\
&& + P_m^{(3)}(x)\frac{2(1-\alpha)}{m-\alpha}x(1-x)
\end{array}
\end{equation}

Let $s_m = \frac{d^2}{dx^2}P_m(x)|_{x=1}$ so that $\sigma_m^2 = s_m + \mu_m - \mu_m^2$.

Evaluating equation (\ref{eqn:cherry_probgenfn_second_derivative}) at $x=1$ gives:
$$\begin{array}{rcl}
s_{m+1} &=& \mu_m\frac{2m(1-\alpha) - 4(1-\alpha)}{m-\alpha} + s_{m}\frac{(m-1)\alpha + m(1-\alpha) - 4(1-\alpha)}{m-\alpha}
\\
&=& \mu_m\frac{2(m-2)(1-\alpha)}{m-\alpha} + s_{m}\frac{m-4 + 3\alpha}{m-\alpha}
\end{array}
$$

Equation (\ref{eqn:cherry_mean_recurrence}) and $\sigma_m^2 = s_m + \mu_m - \mu_m^2$ gives:
$$
\begin{array}{rcl}
s_{m+1} & = & \sigma_{m+1}^2 - \mu_{m+1} + \mu_{m+1}^2 \\
& = & \sigma_{m+1}^2 - \left({ \frac{m(1-\alpha)}{m-\alpha} + \frac{m-2+\alpha}{m-\alpha}\mu_m }\right) + \left({ \frac{m(1-\alpha)}{m-\alpha} + \frac{m-2+\alpha}{m-\alpha}\mu_m }\right)^2 \\
& = & \sigma_{m+1}^2 + \frac{\alpha(1-\alpha)m(1-m)}{(m-\alpha)^2} + \frac{(m-2+\alpha)(m- 2m\alpha) +\alpha)}{(m-\alpha)^2} \mu_m + \frac{(m-2+\alpha)^2}{(m-\alpha)^2}\mu_m^2
\end{array}
$$

Substituting for $s_i$ gives:
\begin{equation}
\label{eqn:sigma_recurrence}
\begin{array}{rcl}
\sigma_{m+1}^2 & = &
\mu_m \frac{2(m-2)(1-\alpha)}{m-\alpha} + \frac{m-4+3\alpha}{m-\alpha}\left( \sigma_m^2 - \mu_m + \mu_m^2\right)
\\ &&
  - \frac{\alpha(1-\alpha)m(1-m)}{(m-\alpha)^2} - \frac{(m-2+\alpha)(m-2m\alpha +\alpha)}{(m-\alpha)^2} \mu_m - \frac{(m-2+\alpha)^2}{(m-\alpha)^2}\mu_m^2
\\
&=&
\frac{\alpha(1-\alpha)m(m-1)}{(m-\alpha)^2} + \sigma_m^2\left(\frac{m-4+3\alpha}{m-\alpha} \right)
\\ &&
+ \mu_m \left(\frac{2(1-\alpha)(m-2\alpha m +\alpha)   }{(m-\alpha)^2}\right)
- \mu_m^2 \left(\frac{4(1-\alpha)^2  }{(m-\alpha)^2}\right)
\end{array}
\end{equation}

From Proposition \ref{prop:cherry-normal}, $\mu_m = \frac{1-\alpha}{3-2\alpha}m + r(m)$ and $\sigma_m^2 = c m + p(m)$, where $r$ and $p$ are $o(m)$ and $c$ is some constant depending on $\alpha$.

Substituting this into equation (\ref{eqn:sigma_recurrence}) and multiplying by $(m-\alpha)^2(3-2\alpha)^2$ gives a quadratic in $m$ which must equal zero.   As $r$ and $p$ are $o(m)$, the coefficient of $m^2$ must tend to zero:
$$c(3-2\alpha)^2(4\alpha-5) - (p(m+1)-p(m))(3-2\alpha)^2 + (1-\alpha)(2-\alpha) \rightarrow 0$$
This gives:
$$\frac{p(m+1)-p(m)}{5-4\alpha} \rightarrow \frac{(1-\alpha)(2-\alpha)}{(3-2\alpha)^2(5-4\alpha)} -c$$
As $c$ is a constant (for fixed $\alpha$) this means that $p(m+1)-p(m)$ must have a limit, which can only by $0$ as $p = o(m)$.
Thus $$c = \frac{(1-\alpha)(2-\alpha)}{(3-2\alpha)^2(5-4\alpha)}$$

\end{proof}

It would of course be nice to know exactly how fast the convergence of $\mu_m$ and $\sigma^2$ is.   More explicit formula are given below.

\begin{corollary}
For $m \ge 3$, $\alpha\in[0,1)$, the expected number of cherries in a random Alpha Tree with $m$ leaves is:
$$\mu_m = \frac{1-\alpha}{3-2\alpha}(m-\alpha) + \frac{\alpha}{2} +
\frac{\alpha}{2(3-2\alpha)}\prod_{i=3}^{m-1} \frac{i-2+\alpha}{i-\alpha}$$
\end{corollary}

\begin{proof}
Let $\mu_m = \frac{1-\alpha}{3-2\alpha}(m-\alpha) + \frac{\alpha}{2} + X_m$.

Then 
$$
\begin{array}{rcl}
\mu_{m+1} &=&
\frac{m(1-\alpha)}{m-\alpha} + \frac{m-2+\alpha}{m-\alpha} \left( \frac{1-\alpha}{3-2\alpha}(m-\alpha) + \frac{\alpha}{2} + X_m \right)
\\ &=&
\frac{1-\alpha}{3-2\alpha}(m+1-\alpha) - (1-\alpha) + \frac{m(1-\alpha)}{m-\alpha} +  \frac{\alpha}{2} \left( 1 - \frac{2(1-\alpha)}{m-\alpha} \right) + \frac{m-2+\alpha}{m-\alpha}X_m
\\ &=&
\frac{1-\alpha}{3-2\alpha}(m+1-\alpha) + \frac{m(1-\alpha) - (m-\alpha)(1-\alpha) - \alpha(1-\alpha) }{m-\alpha} + \frac{\alpha}{2}
+ \frac{m-2+\alpha}{m-\alpha}X_m
\\ &=&
\frac{1-\alpha}{3-2\alpha}(m+1-\alpha) +  \frac{\alpha}{2} + \frac{m-2+\alpha}{m-\alpha}X_m
\end{array}
$$

So $X_{m+1} = \frac{m-2+\alpha}{m-\alpha}X_m$, and so for $m\ge 1$ (and $\alpha \ne 1$)

$$\mu_m = \frac{1-\alpha}{3-2\alpha}(m-\alpha) + \frac{\alpha}{2} +
\prod_{i=1}^{m-1} \frac{i-2+\alpha}{i-\alpha} X_1$$
and $X_1 = 0 - \frac{(1-\alpha)(1-\alpha)}{3-2\alpha} - \frac{\alpha}{2} = \frac{\alpha -2}{2(3-2\alpha)}$.
As $X_3 = \frac{\alpha}{2(3-2\alpha)}$, a more pleasing formula for $m \ge 3$ and all values of $\alpha$ is:

$$\mu_m = \frac{1-\alpha}{3-2\alpha}(m-\alpha) + \frac{\alpha}{2} +
\frac{\alpha}{2(3-2\alpha)}\prod_{i=3}^{m-1} \frac{i-2+\alpha}{i-\alpha}$$

\end{proof}

For rational values of $\alpha$ the product term telescopes.
In the case of $\alpha = 0$, the Yule model,  the expected number of cherries is
$\mu_m = \frac{m}{3}$
In the case of $\alpha = \frac{1}{2}$, the Uniform distribution on cladograms, the expected number of cherries is
$\mu_m = \frac{m(m-1)}{2(2m-3)}$

Note that this second value differs slightly from the numbers given in \cite{McKenzie-Steel-2000} and \cite{Hendy-Penny-1982} as the uniform trees considered there are unrooted.

\section[The shape of evolution: Treebase]{The shape of evolution: Treebase and the big picture}
\label{chapter:treebase}
\subsection{Questions about shape}
This section addresses the shape of phylogenetic trees found in nature and possible biases in common reconstruction techniques.
The recent increase in protein and nucleotide sequence data and availability of programs for reconstructing phylogeny from such data has lead to a large number of published phylogeny.  Many of these phylogenetic trees have been made available in online databases, such as Treebase \cite{Treebase}.

Some natural questions that arise are: How asymmetrical are the trees found in nature?  Do they follow some nice probability distribution and if so what is it?  Are all trees about the same shape?  Are there systematic biases in different reconstruction techniques?  An excellent discussion of these issues is given by Mooers and Heard \cite{Mooers-Heard-1997}.

The question of the `amount of asymmetry' in natural trees is often raised.
One major stumbling block in a systematic analysis of tree shapes has been the absence of a good measure of imbalance.  Heard's analysis \cite{Heard-1992} of 208 published phylogeny is hampered by exactly this problem.
Several measures of tree imbalance have been considered in the past such as ``Colless's I'' and ``Sackin's index'' (  see Section \ref{chapter:Sackins-Colless-index} for a description of these statistics).  Unfortunately these statistics change greatly with the number of leaves, and have means and variances depending on the probability distribution chosen (see \cite{Rogers-1996} for example).

It has often been observed that phylogenetic trees found in nature are in general more symmetric than Uniform trees but not as symmetric as Yule trees (for example \cite{Mooers-Heard-1997},\cite{Heard-1996}).  This observation is verified and quantified here by examining the distribution over the trees in Treebase of the maximum likelihood estimate of the parameter in the alpha model.  The median of these estimates is about $\alpha=0.22$, directly between the Yule ($\alpha=0$) and Uniform ($\alpha=0.5$) models.

A variety of statistics are used to measure how close the data fits the alpha model.  Combined p-values are used to reject the hypothesis that the trees in Treebase all fit with the alpha model.

Other than, perhaps, \cite{Heard-1992} this appears to be the first systematic analysis of the shape and balance of a large number of published phylogeny.
\ignore{I indend to continue this analysis in more detail, and include comparisons with fitting the beta model of Aldous.}

\subsection{Estimating alpha}
The probability of a given tree shape under the alpha model is a rational function of alpha, and may be easily computed.  By Lemma \ref{lemma:alpha-mss}, the alpha model is Markovian self-similar with conditional split distribution
$$q_\alpha(a,b) = \frac{\Gamma_\alpha(a)\Gamma_\alpha(b)}{\Gamma_\alpha(a+b)}
\left({ \frac{\alpha}{2}\binom{a+b}{a} + (1-2\alpha)\binom{a+b-2}{a-1} }\right)$$
where $\Gamma_\alpha(n) = (n-1-\alpha)(n-2-\alpha)\cdots (2-\alpha)(1-\alpha)$ and $\Gamma_\alpha(1) = 1$

By Proposition \ref{prop:tree-prob-thin-shapes} the probability of a tree shape under such a model is the product of the conditional split probability at each branch-point.  See Section \ref{sec:probability-of-trees} for more details and an example.

In the analysis presented here the probability of the tree shape was calculated for $1000$ equally spaced values of $\alpha$ in $[0,1]$.  The maximum over these $1000$ points was then taken as a good approximation of the maximum likelihood for alpha.

Two transformations to the set of trees were made before estimating alpha.  The first was to remove all non-binary trees as there they are not covered by the model (and probably indicate insufficient data to reconstruct a tree \cite{Heard-1992}).  The second was to accomodate the fact that most published trees contain an outgroup.

In many phylogenetic reconstructions, an outgroup is used to locate the root on a reconstructed tree, as many algorithms give unrooted trees or unsure root positions.  An outgroup is a singleton or pair (or more) of taxa which are artificially chosen to be significantly different from the rest of the taxa in the analysis.  These are then used to root the reconstructed tree as it is assumed that the first speciation event separates the outgroup taxa from the from the main group, sometimes called the ingroup.

The addition of outgroups in this manner is expected to increase the average imbalance of trees and the maximum likelihood estimate for alpha.  To avoid this bias, all trees were split at the root into two separate trees.  Trees of size 3 or less were all discarded.  In the event that a tree was constructed without an outgroup this should not greatly effect the estimate of alpha, particularly if the tree shape obeys a Markovian self-similar model (as seems evolutionarily plausible).  Almost all trees in the sample set appeared to have an outgroup.

The median values for the maximum likelihood estimates for alpha before and after this splitting at the root were about $0.37$ and $0.22$ respectively.  Thus, removal of the outgroup does significantly effect the estimation of alpha.  This is to be expected for trees of the size most present in Treebase.  Estimation for larger trees should be less effected by the presence of an outgroup.

\begin{figure}[hp]
\center{\includegraphics[width=5in]{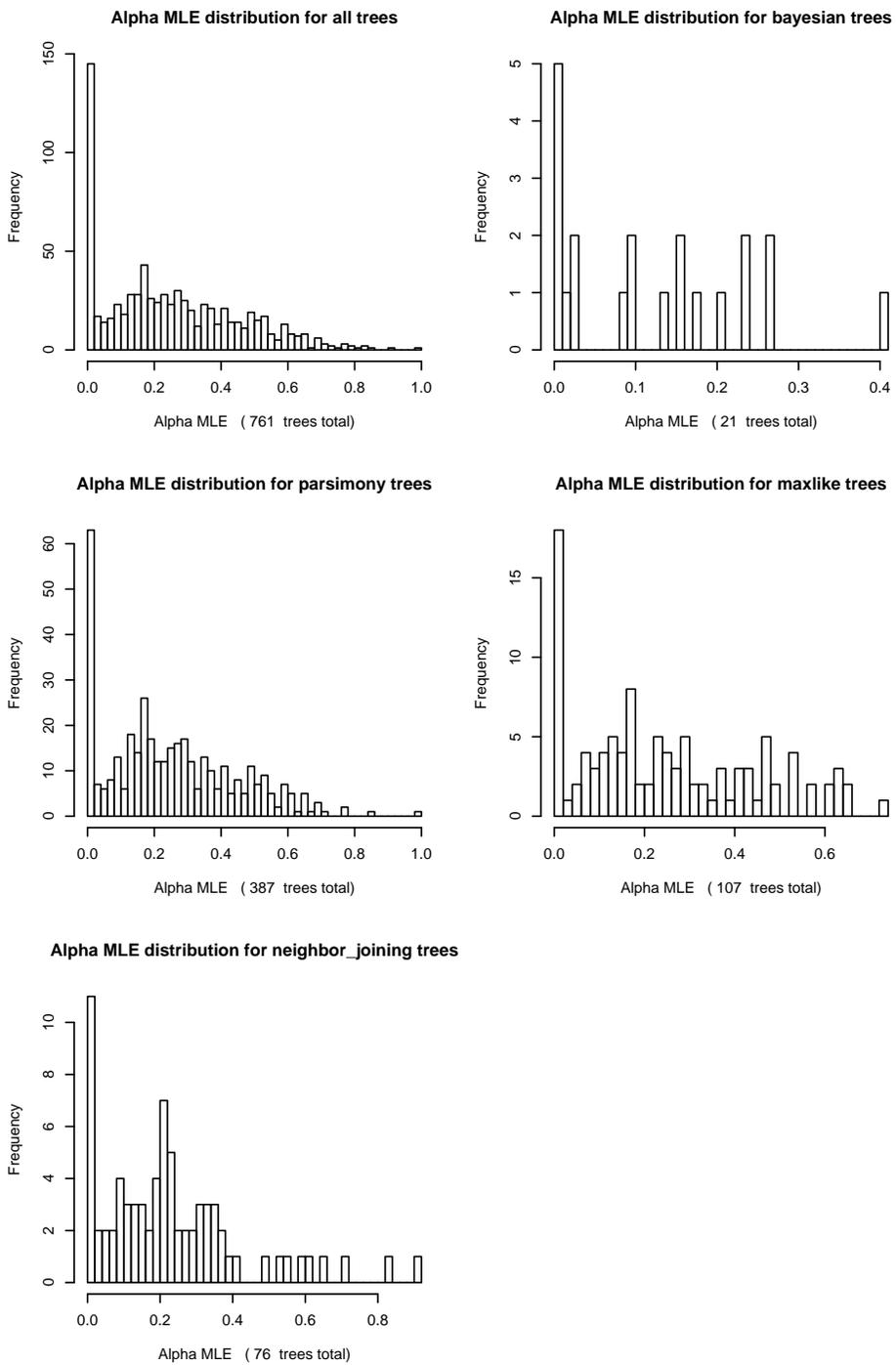}}
\caption{Maximum likelihood estimates of alpha for trees with at least 10 leaves (outgroups removed)}
\label{fig1}
\end{figure}

Figure \ref{fig1} shows a histograph for the maximum likelihood estimates for alpha, categorized by reconstruction method.  All trees with less than 10 leaves were discarded as for small trees the number of different shapes is too small to allow for a fine estimate of alpha.  The number of trees remaining was 761.

Treebase entries also include the method of reconstruction in most cases.
Here are summary statistics, with a break-down by reconstruction method.

\noindent
\begin{tabular}{llllllll}
Method              & \# trees   &   Min. & 1st Qu. & Median &   Mean & 3rd Qu.& Max.\\ 
all & 761               & 0.0000 & 0.0800 & 0.2200 & 0.2536 & 0.3900 & 0.9900 \\
parsimony& 387          & 0.0000 & 0.1050 & 0.2300 & 0.2565 & 0.3800 & 0.9900 \\
maxlike& 107            & 0.0000 & 0.1000 & 0.2300 & 0.2545 & 0.4100 & 0.7400 \\
neighbor joining& 76    & 0.0000 & 0.0975 & 0.2150 & 0.2361 & 0.3225 & 0.9200 \\
bayesian& 21            & 0.0000 & 0.0200 & 0.1000 & 0.1262 & 0.2100 & 0.4100 \\
unknown& 170
\end{tabular}
\break

Note that the median is consistently around $0.22$ (except for the bayesian method).
The number of trees with maximum likelihood estimate for alpha strictly between $0$ and $0.5$ is $511$ out of a possible $761$ (about $67\%$).

Applying a t-test to the estimates for parsimony and bayesian methods gives a p-value of $5.943 *10^{-5}$ (degrees of freedom=26.494).  This indicates a strong differential bias between these two reconstruction methods.
However, it should be noted not all methods were applied to all data.  It may be that phylogenists working on different types of organism with different average tree shapes may prefer one reconstruction method over the other.
In order to do a fully systematic study each method should be applied to the original sequence data where it is available.

The large spike at alpha $=0$ (about $20\%$ of the trees) is discussed in the next section.

\subsection{Does the data fit the model}
This section covers the question of how well the data fits the model.  This is addressed using p-value data for a number of different statistics on trees.
As an explicit model is being tested there is no longer a problem with using statistics which change with the number of leaves.

Given a statistic, model and tree, the distribution of the statistics under the model can be compared with the statistic on the given tree, to give a p-value.  If the trees are generated by the model then such p-values are uniformly distributed (at least for continuous distributions).

For each tree, and statistic, this p-value was estimated by generating $1000$ random trees from the model (with the MLE value of alpha) to approximate the distribution under the model.
This estimate has the correct mean, and a variance of at most $\frac{1}{4\sqrt{1000}}$.

The statistics used are "Colless' I", "number of cherries" (pairs of adjacent leaves), "total depth of all leaves" (Sackin's Index) (equivalently: average leaf depth), "maximum depth of a leaf", and the probability (considered as a function on the set of trees of a fixed size).
See Sections \ref{chapter:Sackins-Colless-index} and \ref{chapter:cherries} for more details on these statistics.

\begin{figure}[hp]
\center{\includegraphics[width=5in]{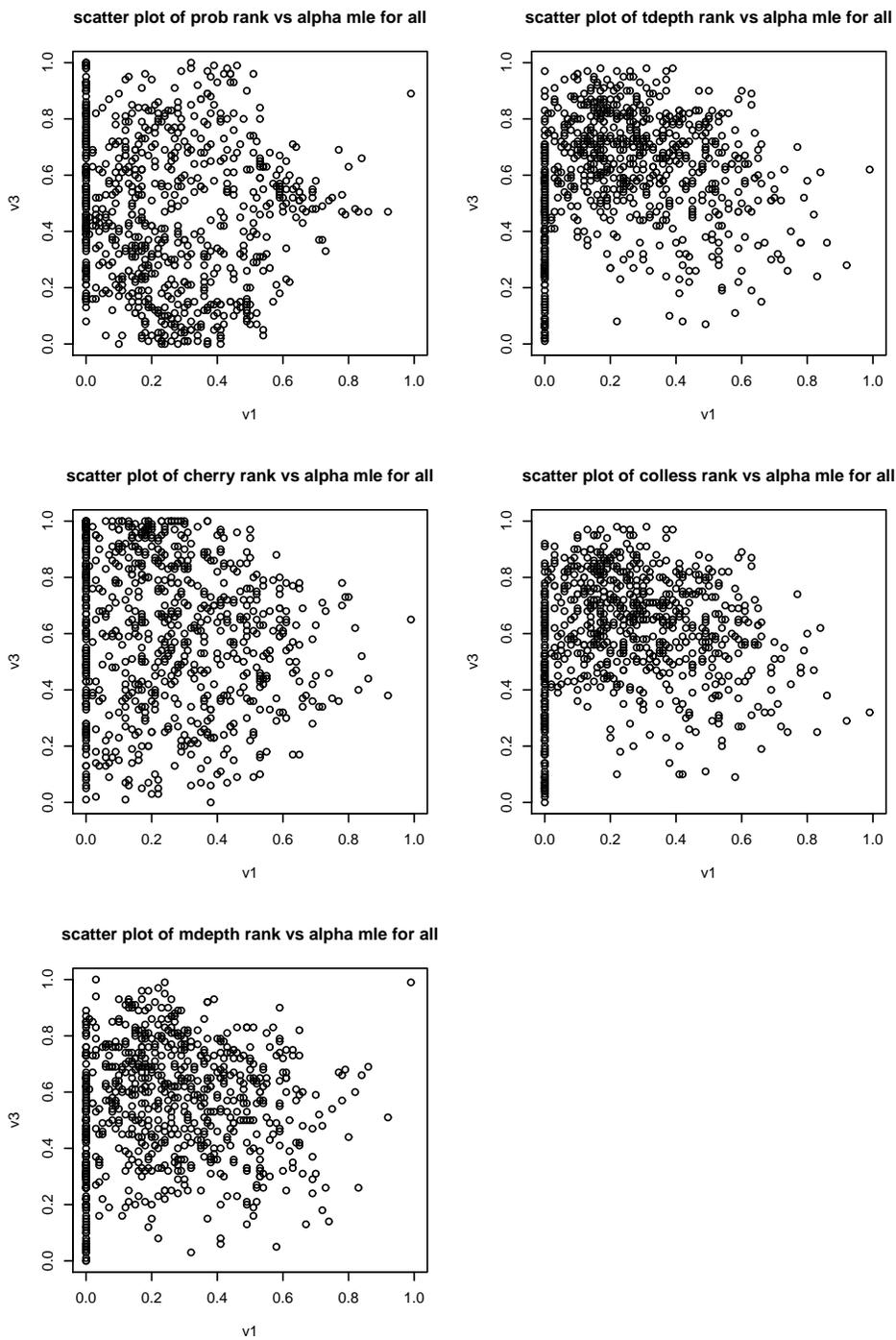}}
\caption{P-values for various statistics plotted against alpha MLE}
\label{fig2}
\end{figure}

\begin{figure}[hp]
\center{\includegraphics[width=5in]{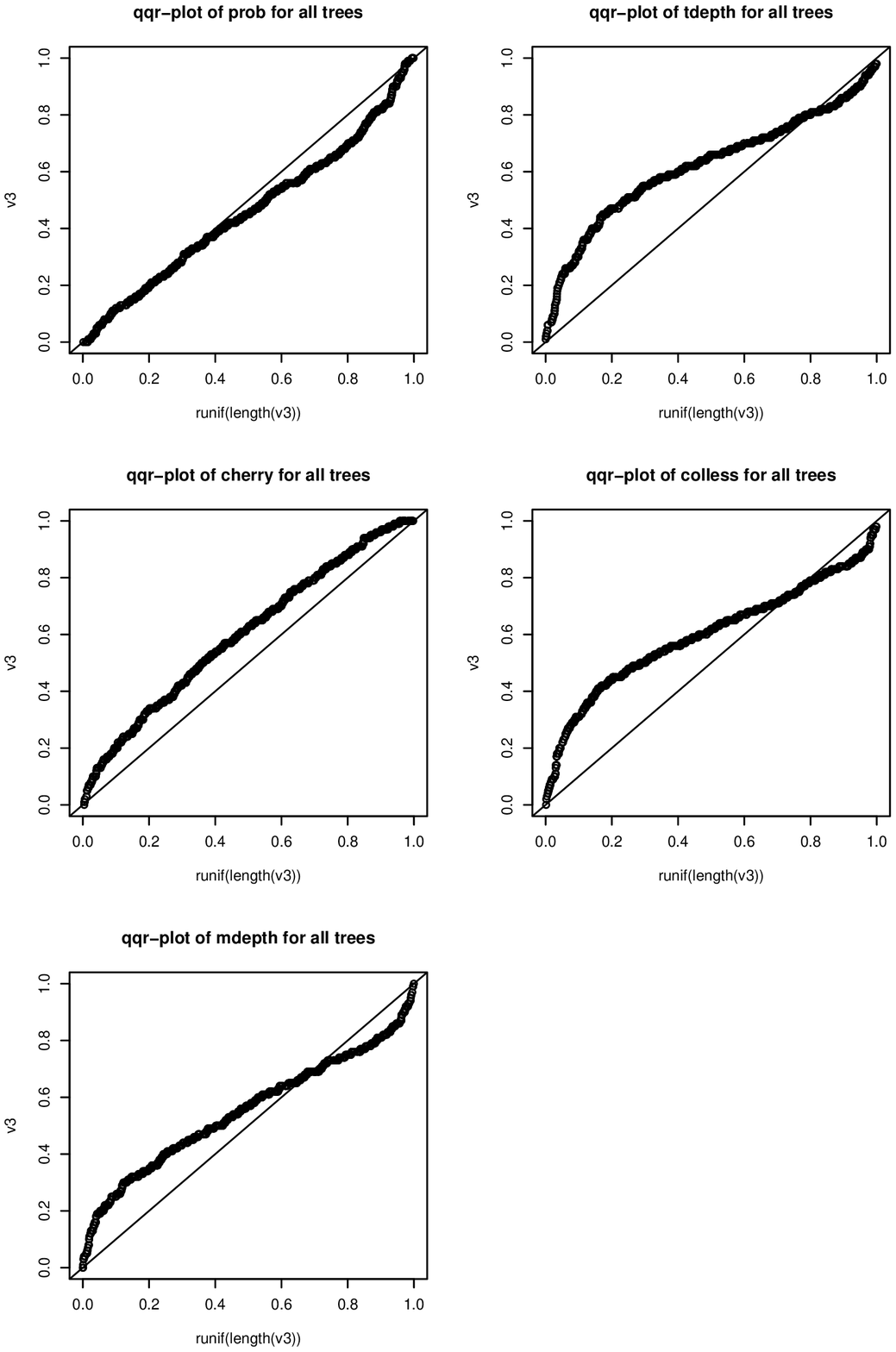}}
\caption{Q-Q Plots of for various statistics}
\label{fig3}
\end{figure}

Figure \ref{fig2} shows scatter plots of these p-values against the estimate of alpha.  Figure \ref{fig3} shows qq-plots of these p-values against uniform $[0,1]$.

Looking at these plots it is clear that while the model is not terrible, it is certainly not a perfect fit.  The lack of extreme p-values for alpha in $(0,1)$ might be explained by extreme trees being better fit by other values of alpha where their shape is not so unusual or extreme.  This may also explain some of the large spike at $\alpha=0$, which comprises about $20\%$ of all the trees.

\clearpage
\section{Tree shapes with up to $7$ leaves}
\label{chapter:tree-shapes-up-to-7}
This appendix contains a list of all tree shapes with up to $7$ leaves, ordered lexicographically.
The number of phylogenetic trees with a given shape is stated, as well as the probability of this shape under the alpha model (conditional on the tree having that many leaves).

The probability of a tree shape $T$ under the alpha model is, by Proposition \ref{prop:tree-prob-thin-shapes}, equal to:
$$ \prod_{(a,b)\in \{\text{splits}(T)} \hat q_\alpha\{a,b\}$$
Recall that $\hat q_\alpha\{a,b\} = q_\alpha(a,b) + q_\alpha(b,a)$ if $a\ne b$ and $\hat q_\alpha\{a,a\} = q_\alpha(a,a)$.
Equation \ref{eqn:alpha-split-distribution} provides the split distribution, $q_\alpha$, of the alpha model.

If $A(n)$ is the number of tree shapes with exactly $n$ leaves then $A(m)$ satisfies the following recurrence relations: $A(2n+1) = \sum_{i=1}^{n} A(i)A(2n+1-i)$, $A(2n) = \sum_{i=1}^{n-1} A(i)A(2n-i) + \frac{A(n)(A(n)-1)}{2} + A(n)$.  Set $A(0)=0$ for convenience and $A(1)=1$.
Thus, the first few values of $A(n)$ are  1, 1, 1, 1, 2, 3, 6, 11, 23, 46, 98, 207, 451, 983, 2179, 4850, 10905, 24631, 56011, 127912, 293547, 676157.

This is sequence A001190 in the Encyclopedia of integer sequences.  The generating function, $G(x)$, of sequence A001190 satisfies the recurrence relation $G(x) = x + (1/2)(G(x)^2 + G(x^2))$

\epsfysize=2 cm
\begin{figure}[htbp]
\center{
	\includegraphics{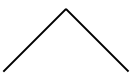}
	\caption[Treeshape $(1,1)$]{The trivial two-leaf tree, treeshape $(1,1)$.
	 \newline There is $1$ phylogenetic tree with this shape.
	 \newline The probability of this shape under the alpha model is $1$}
}
\end{figure}

\epsfysize=2 cm
\begin{figure}[htbp]
\center{
	\includegraphics{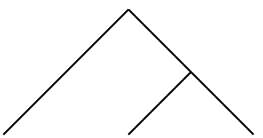}
	\caption[Treeshape $(2,1)$]{The unique three-leafshape, treeshape $(2,1)$
	 \newline There are $\frac{3!}{2}=3$ phylogenetic tree with this shape.
	 \newline The probability of this shape under the alpha model is $1$}
}
\end{figure}

\epsfysize=2 cm
\begin{figure}[htbp]
\center{
	\leavevmode
	\epsfbox{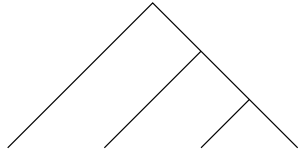}
	\caption[Treeshape $(3,1)$]{The four-leaf comb, treeshape $(3,1)$
	 \newline There are $\frac{4!}{2} = 12$ phylogenetic trees with this shape. 	 \newline The probability of this shape under the alpha model is $\frac{2}{3-\alpha}$}
}
\end{figure}

\epsfysize=2 cm
\begin{figure}[htbp]
\center{
	\leavevmode
	\epsfbox{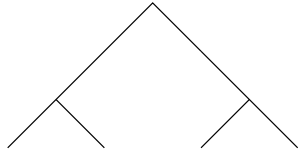}
	\caption[Treeshape $(3,2)$]{Treeshape $(3,2)$, sequence $4211211$ or $4$.
	 \newline There are $\frac{4!}{2^3} = 3$ phylogenetic trees with this shape.
	 \newline The probability of this shape under the alpha model is $\frac{1-\alpha}{3-\alpha}$}
}
\end{figure}

\epsfysize=2 cm
\begin{figure}[htbp]
\center{
	\leavevmode
	\epsfbox{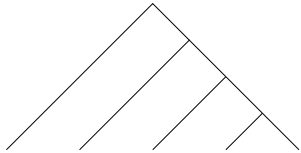}
	\caption[Treeshape $(4,1)$]{Treeshape $(4,1)$, sequence $543211111$ or $543$
	\newline There are $\frac{5!}{2} = 60$ phylogenetic trees with this shape.
	 \newline The probability of this shape under the alpha model is $\frac{2(2+\alpha)}{(4-\alpha)(3-\alpha)}$}
}
\end{figure}

\epsfysize=2 cm
\begin{figure}[htbp]
\center{
	\leavevmode
	\epsfbox{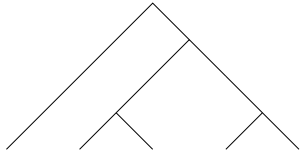}
	\caption[Treeshape $(4,2)$]{Treeshape $(4,2)$, sequence $542112111$ or $54$.
	 \newline There are $\frac{5!}{2^3} = 15$ phylogenetic trees with this shape.
	 \newline The probability of this shape under the alpha model is $\frac{(1-\alpha)(2+\alpha)}{(4-\alpha)(3-\alpha)}$}
}
\end{figure}

\epsfysize=2 cm			
\begin{figure}[htbp]
\center{
	\leavevmode
	\epsfbox{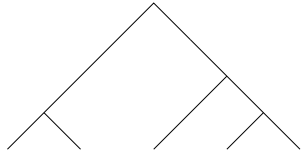}
	\caption[Treeshape $(4,3)$]{Treeshape $(4,3)$, sequence $532111211$ or $53$.
	\newline There are $\frac{5!}{2^2} = 30$ phylogenetic trees with this shape.
	 \newline The probability of this shape under the alpha model is $\frac{2(1-\alpha)}{4-\alpha}$}
}
\end{figure}

\epsfysize=2 cm			
\begin{figure}[htbp]
\center{
	\leavevmode
	\epsfbox{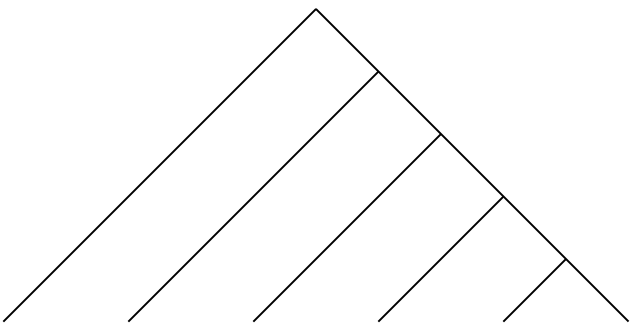}
	\caption[Treeshape $(5,1)$]{Treeshape $(5,1)$, sequence $65432111111$ or $6543$.
	 \newline There are $\frac{6!}{2} = 360$ phylogenetic trees with this shape.
	 \newline The probability of this shape under the alpha model is
	 $\frac{4(1+\alpha)(2+\alpha)}{(5-\alpha)(4-\alpha)(3-\alpha)}$}
}
\end{figure}

\epsfysize=2 cm			
\begin{figure}[htbp]
\center{
	\leavevmode
	\epsfbox{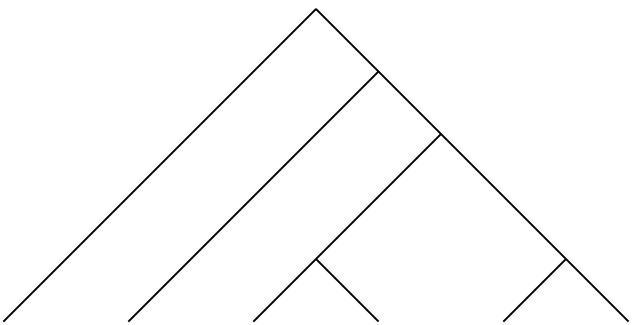}
	\caption[Treeshape $(5,2)$]{Treeshape $(5,2)$, sequence $65421121111$ or $654$.
	 \newline There are $\frac{6!}{2^3} = 90$ phylogenetic trees with this shape.
	 \newline The probability of this shape under the alpha model is 
	 $\frac{2(1-\alpha)(1+\alpha)(2+\alpha)}{(5-\alpha)(4-\alpha)(3-\alpha)}$}
}
\end{figure}

\epsfysize=2 cm			
\begin{figure}[htbp]
\center{
	\leavevmode
	\epsfbox{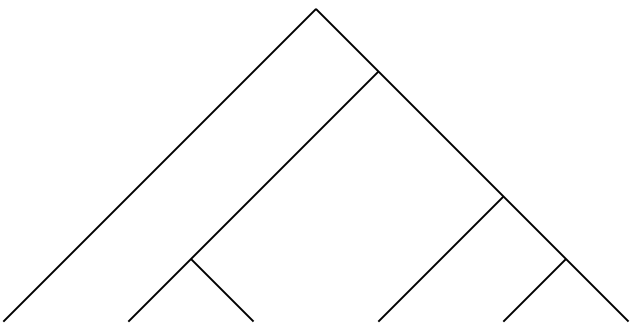}
	\caption[Treeshape $(5,3)$]{Treeshape $(5,3)$, sequence $6532111211$ or $653$.
	\newline There are $\frac{6!}{2^2} = 180$ phylogenetic trees with this shape.
	\newline The probability of this shape under the alpha model is
	$\frac{4(1-\alpha)(1+\alpha)}{(5-\alpha)4-\alpha)}$}
}
\end{figure}

\epsfysize=2 cm			
\begin{figure}[htbp]
\center{
	\leavevmode
	\epsfbox{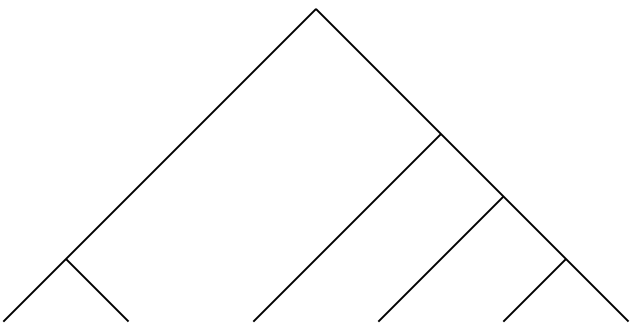}
	\caption[Treeshape $(5,4)$]{Treeshape $(5,4)$, sequence $64321111211$ or $643$.
	  \newline There are $\frac{6!}{2^2} = 180$ phylogenetic trees with this shape.
	\newline The probability of this shape under the alpha model is 
	$\frac{2(1-\alpha)(8-\alpha)}{(5-\alpha)(4-\alpha)(3-\alpha)}$}
}
\end{figure}

\epsfysize=2 cm			
\begin{figure}[htbp]
\center{
	\leavevmode
	\epsfbox{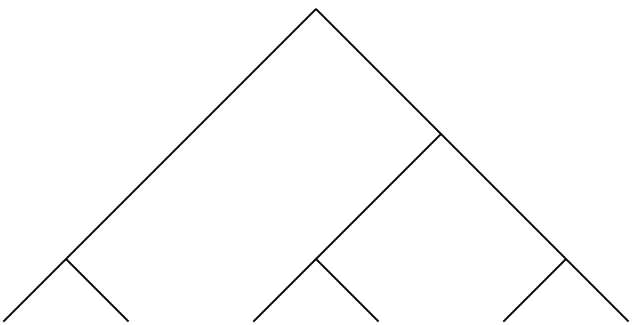}
	\caption[Treeshape $(5,5)$]{Treeshape $(5,5)$, sequence $64211211211$ or $64$.
	  \newline There are $\frac{6!}{2^4} = 45$ phylogenetic trees with this shape.
	\newline The probability of this shape under the alpha model is 
	$\frac{(1-\alpha)^2(8-\alpha)}{(5-\alpha)(4-\alpha)(3-\alpha)}$}
}
\end{figure}

\epsfysize=2 cm			
\begin{figure}[htbp]
\center{
	\leavevmode
	\epsfbox{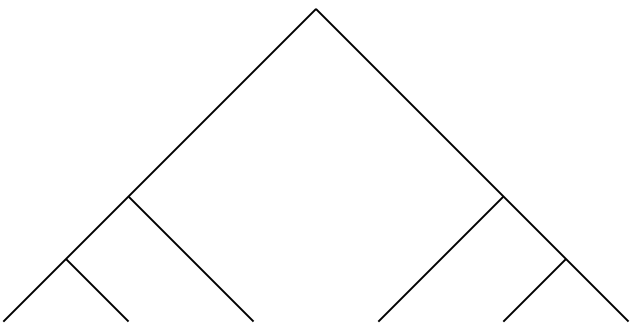}
	\caption[Treeshape $(5,6)$]{Treeshape $(5,6)$, sequence $63211132111$ or $633$.
	  \newline There are $\frac{6!}{2^3} = 90$ phylogenetic trees with this shape.
	\newline The probability of this shape under the alpha model is 
	$\frac{2(2-\alpha)(1-\alpha)}{(5-\alpha)(4-\alpha)}$}
}
\end{figure}

\bibliography{phylo}
\bibliographystyle{plain}
\end{document}